\documentclass[oupthm,crop,info]{CUP-JNL-CJM}%
\usepackage{amsmath, amssymb, amsthm, mathrsfs, enumerate, tikz, BOONDOX-cal, units, enumitem, hyperref}
\usetikzlibrary{positioning, calc, matrix, arrows.meta}
\usepackage[utf8]{inputenc}

\theoremstyle{oupplain}
\newtheorem{theorem}{Theorem}[section]
\newtheorem{prop}[theorem]{Proposition}
\newtheorem{cor}[theorem]{Corollary}
\newtheorem{lem}[theorem]{Lemma}
\theoremstyle{oupdefinition}
\newtheorem{definition}{De{f}inition}[section]
\theoremstyle{oupremark}
\newtheorem{remark}[theorem]{Remark}
\newtheorem{example}[theorem]{Example}
\theoremstyle{oupproof}
\newtheorem{proof}{Proof}

\numberwithin{equation}{section}

\articletype{RESEARCH ARTICLE}
\jname{Canadian Mathematical Society}
\jyear{2021}

\def\C*{{\sl C*}-algebra} 
\def\Cs*{{\sl C*}-subalgebra} 
\def\HC*{{Hilbert {\sl C*}}} 

\newcommand{\CC}{{\mathbb{C}}}
\newcommand{\NN}{{\mathbb{N}}}
\newcommand{\norm}[1]{\lVert#1\rVert} 
\newcommand{\bignorm}[1]{\bigl\lVert#1\bigr\rVert} 
\newcommand{\setst}[2]{\ensuremath{\left\{{#1}\,|\,{#2} \right\}}} 
\newcommand{\Hom}[2]{\ensuremath\textup{Hom}{({#1},{#2})}}        
\newcommand{\Homunder}[3]{\ensuremath\textup{Hom}_{#3}{({#1},{#2})}} 
\newcommand{\Mor}[2]{\ensuremath\textup{Mor}{({#1},{#2})}}
\newcommand{\Morunder}[3]{\ensuremath\textup{Mor}_{#3}{({#1},{#2})}}
\newcommand{\cats}[1]{\ensuremath{\mathcal{#1}}} 
\newcommand{\catsa}{\cats{A}}

\newcommand{\mapfromto}[3]{{#1}\colon {#2}\rightarrow {#3}}

\newcommand{\cbnorm}[1]{\norm{#1}_{\textup{cb}}}
\newcommand{\idntyof}[1]{\ensuremath{\textup{id}}_{#1}} 

\newcommand\CB[2]{C\!B({#1},\,{#2})}

\newcommand\CBmod[3]{C\!B_{#1}({#2},\,{#3})}

\newcommand\classofmorphisms[3]{{\mathcal{#1}}{({#2},\,{#3})}}
\newcommand\morphismsnomathcal[3]{{{#1}}{({#2},\,{#3})}}
\newcommand\classMofmorphisms[2]{{\mathcal{M}}{({#1},\,{#2})}}

\newcommand{\Abcat}{{\mathcal{A\!b}}}   
\newcommand{\Ex}{{\mathcal{E\!x}}}   
\newcommand{\Exmax}{{\Ex_{max}}}
\newcommand{\Exmin}{{\Ex_{min}}}
\newcommand{\Exrel}{{\Ex_{rel}}}
\newcommand{\Mrel}{{\mathcal{M}_{rel}}}
\newcommand{\Nrel}{{\mathcal{N}_{rel}}}
\newcommand{\Prel}{{\mathcal{P}_{rel}}}

\newcommand{\Banone}{{\mathcal{B\!a\!n}^{\lower2pt\hbox{\tiny$1$}}}}
\newcommand{\Baninf}{{\mathcal{B\!a\!n}^{\lower1pt\hbox{\tiny$\infty$}}}}
\newcommand{\BanAinf}{{\mathcal{B\!a\!n}^{\lower1pt\hbox{\tiny$\infty$}}_A}}
\newcommand{\Opone}{{\mathcal{O\!p}^{\lower1.5pt\hbox{\tiny$1$}}}}   
\newcommand{\Opinf}{{\mathcal{O\!p}^{\lower.5pt\hbox{\tiny$\infty$}}}}  

\newcommand{\OMODinf}[1]{{\ensuremath{\mathcal{O\!M\!o\!d}^\infty_{{#1}}}}}

\newcommand{\OMODAinf}{{\OMODinf{\!A}}}
\newcommand{\OMODCCinf}{{\OMODinf{\!\CC}}}

\newcommand{\XMODinfA}{{\ensuremath{\mathcal{X\!M\!o\!d}^\infty_{{A}}}}}

\newcommand{\Mmor}{{\mathcal{M}}}  

\newcommand{\mnMODinf}[1]{{\ensuremath{\mathcal{m\!n\!M\!o\!d}^\infty_{{#1}}}}}

\newcommand{\mnMODAinf}{{\mnMODinf{\!A}}}

\newcommand{\Ffunct}{{\mathsf{F}}}

\newcommand{\tenmod}[3]{{#1}\otimes_{{#2}}\mkern2mu{#3}}
\DeclareMathOperator*{\otimesfrown}{\!{}\otimes \mkern-15.8mu{{\text{\raisebox{1.35ex}{\tiny{$\frown$}}}}}{}}
\DeclareMathOperator*{\otimesfrowna}{\!{}\otimes \mkern-15.6mu{{\text{\raisebox{1.35ex}{\tiny{$\frown$}}}}}{}}
\newcommand{\projten}[2]{{#1}\otimesfrown\mkern2mu{#2}}   
\newcommand{\projtenmod}[3]{{#1}\otimesfrowna\nolimits_{{#2}}\mkern2mu{#3}}    
\newcommand{\hten}[2]{\tenmod{{#1}}{h}{{#2}}} 
\newcommand{\htenmod}[3]{\tenmod{{#1}}{h{#2}}{{#3}}} 

\newcommand{\Midim}[1]{\ensuremath\textup{Inj}_{\mathcal{M}}\textup{-dim}\,{(#1)}} 
\newcommand{\dgcb}[1]{{\ensuremath\textup{dg}_{cb}\,{#1}}}  

\let\oldmarginpar\marginpar
\renewcommand\marginpar[1]{\-\oldmarginpar[\raggedleft\bf\scriptsize #1]%
{\raggedright\bf\scriptsize #1}}

\begin{document}

\begin{Frontmatter}

\title{Exact Structures for Operator Modules}
\author{Martin Mathieu}
\author{Michael Rosbotham}
\address{\orgname{Mathematical Sciences Research Centre, Queen's University Belfast},
\orgaddress{Belfast BT7 1NN}, \country{Northern Ireland},
\email{m.m@qub.ac.uk},
\email{mrosbotham01@qub.ac.uk}}


\abstract{%
We demonstrate how exact structures can be placed on the additive category of right operator modules over an operator algebra
in order to discuss global dimension for operator algebras. The properties of the Haagerup tensor product play a decisive role
in this.}

\end{Frontmatter}

\maketitle

\section{Introduction}\label{sect:intro}

\noindent
Among the most important operator space modules over \C*s are the \HC*-modules, the operator modules and the matrix normed modules.
The first class became prominent through the work of Paschke and Rieffel and was intensively studied by Blecher, amongst others,
see, e.g., \cite{Blecher1997} and~\cite{BLM04}. The main difference between the second and the third class lies in the kind
of complete boundedness which is required of the bilinear mappings that give the module action. For operator modules
(we will follow the terminology of~\cite{BLM04} in this paper), one demands (multiplicatively) completely bounded bilinear mappings,
and the associated tensor product is the Haagerup module tensor product~$\htenmod{}{A}{}$ over the \C*~$A$.
Jointly completely bounded bilinear mappings and the module operator space projective tensor product~$\projtenmod{}{A}{}$
govern the class of matrix normed modules; for details, we refer to \cite[Chapter~3]{BLM04}.
Both these classes have been put to good use and found a range of interesting applications; we only mention the recent papers
\cite{Bearden-Crann}, \cite{Crann2017} and \cite{Crisp2018}, \cite{CrispHigson2016} as samples.

For historical reasons the terminology is (still) not uniform; we will follow the convention in~\cite{BLM04}, see also
\cite[Appendix~A]{Rosb2020}. Suffice it to say at this point that operator modules form a full subcategory of the matrix normed modules
(over any operator algebra). These categories are not abelian and therefore the usual homological algebra does not apply directly.
Nevertheless, homological methods have been successfully developed in this framework, for example by Helemskii and his school.
In \cite[Theorem III.5.17, Corollary IV.4.12]{HELhomBook} Helemskii proved that a \C* $A$ is classically semisimple
(i.e., finite dimensional) if and only if all objects in the category $\BanAinf$ of Banach $A$-modules are projective
with respect to the class of epimorphisms that split as morphisms in~$\Baninf$. See also \cite{Aristov2000}, \cite{Aristov2002}
and~\cite{Helem2009}.

For an operator algebra $A$ (on a Hilbert space) it is pertinent to use operator modules to build a cohomology theory
(for the definitions, see Section~\ref{Subsect:op-spaces}); of the numerous contributions we only mention~\cite{BlMuPa2000}, \cite{Helem2009},
\cite{Paulsen1998}, \cite{Volo2009} and~\cite{Wood1999} here.
In this paper, we focus on an appropriate definition of cohomological dimension and, in particular, answer a question raised
by Helemskii~\cite{Helem2009} whether quantised global dimension zero is equivalent to the algebra being classically semisimple;
see Theorem~\ref{Theorem: Relative Global Dimension 0 iff classically semisimple} below.
In contrast to the situation in ring theory, it appears necessary to limit ourselves to a relative cohomology theory since,
otherwise, there exist too many monomorphisms (equivalently, epimorphisms) and the concepts of injectivity (respectively, projectivity)
become too restrictive. For example, the canonically defined injective global dimension of any \C*
with regards to the category of its operator modules always is at least~$2$, as shown in~\cite{Rosb2020}.
Paulsen discussed relative cohomology in~\cite{Paulsen1998} and related it to completely bounded Hochschild cohomology.
His notions of relative injectivity and projectivity coincide with ours (defined in
Section~\ref{Section: Exact Categories and Global Dimension} below) and Helemskii's approach is also compatible.

The main novelty in this paper is the systematic use of exact categories enabling us to bring category theory to the foreground
in order to take advantage of its unifying strength. Typically, categories that appear in analysis are not abelian which has hindered
the full application of homological algebra methods. Exact structures in the sense of Quillen~\cite{Quill73} (to be defined in
Section~\ref{Section: Exact Categories and Global Dimension}) on additive categories like $\Baninf$ were employed by B\"uhler in~\cite{B11},
see also~\cite{B08} for a very nice introduction.
The use of exact structures for categories of operator modules, and indeed sheaves of operator modules over \textsl{C*}-ringed spaces,
was initiated in~\cite{AM2}, see also the survey articles~\cite{M19} and~\cite{MR20}, and further studied in~\cite{Rosb2020}.
Exact categories are ideally suited for the important tools of homological algebra such as (short) exact sequences, diagram
lemmas and derived functors.

Not assuming an in-depth knowledge of category theory (as this article is written for functional analysts)
we have included the main categorical concepts deployed throughout this paper in Section~\ref{Subsection: Category Theory}.
In Section~\ref{Subsect:op-spaces} our assumptions on operator spaces, operator algebras and the type of operator space modules we work with
are stated.

In Section~\ref{Section: Additive category of operator modules} we show how the categorical notions of kernels and cokernels can be expressed,
using language and concepts of operator space theory, in the additive category $\OMODAinf$ of non-degenerate right operator $A$-modules
over an operator algebra~$A$. This prepares Section~\ref{Section: Exact Categories and Global Dimension} in which
we explore the concept of exact categories and some techniques for working in this setting.
We show that there is a canonical exact structure that can be placed on our main category $\OMODAinf$
(Theorem~\ref{Theorem: Emax in OMODAinf forms an exact structure}) and discuss how a global dimension arises in exact categories
with enough injectives.

Section~\ref{Section: Relative Homological Algebra for Operator Modules over Operator Algebras} deals with another exact structure
that can be placed on $\OMODAinf$. This structure is related to the study of `relative homological algebra' for operator modules
as it is done in~\cite{Aristov2002}, \cite{Helem2009} and~\cite{Paulsen1998}, for example.
Our main Theorem~\ref{Theorem: Relative Global Dimension 0 iff classically semisimple} states that, for a unital operator algebra $A$,
the global dimension of this exact category is zero if and only if $A$ is classically semisimple, hence a finite direct sum
of full matrix algebras over the complex numbers.
We also discuss similarities and differences between our preferred category, $\OMODAinf$, and $\mnMODAinf$, the category of matrix normed
modules over an operator algebra~$A$ in this section.

\section{Terminology, Notation and Conventions}\label{Section: Terminology}
\subsection{Category Theory}\label{Subsection: Category Theory}
Let $\mathcal{A}$ be a category. To indicate that $E$ belongs to the class of objects in $\mathcal{A},$ we write $E\in\mathcal{A}$.
The set of morphisms from $E\in\mathcal{A}$ to $F\in\mathcal{A}$ is denoted $\Morunder{E}{F}{\mathcal{A}}$,
with the $\mathcal{A}$ sometimes dropped if the category we are working in is obvious.
If $\mathcal{M}$ is a class of morphisms in $\mathcal{A}$ we denote the subset of $\Morunder{E}{F}{\mathcal{A}}$ that
consists only of the morphisms in $\mathcal{M}$ by $\classMofmorphisms{E}{F}$. Two of the most important classes
of morphisms in a category are the classes of monomorphisms and of epimorphisms. Recall that a morphism $f$ in $\mathcal{A}$
is a \textit{monomorphism} if it is `left cancellable', that is, if $g,h$ are morphisms in $\mathcal{A}$, composable with $f$,
such that $fg=fh$, then we must have $g=h$.  The `right cancellable' morphisms in a category are \textit{epimorphisms}.
An important class of epimorphisms are the retractions. A morphism $r\in\Morunder{E}{F}{\mathcal{A}}$ is a \textit{retraction}
if there exists a morphism $s\in\Morunder{F}{E}{\mathcal{A}}$ such that $rs=\idntyof{F}$, the identity morphism of~$F$.
In this case we say $s$ is a \textit{section} of $r$ and $F$ is a \textit{retract} of~$E$.

The categories we work in will contain a zero object and, therefore, zero morphisms. Thus we can talk about kernels and cokernels.

\begin{definition}\label{Definition of a kernel and a cokernel}
Let $\mathcal{A}$ be a category with a zero object. Suppose $E,F\in\cats{A}$ and  $f\in\Mor{E}{F}$.
	
A \textit{kernel} of $f$ is a pair $(K,\mu)$, where $K\in\cats{A}$ and $\mu\in\Mor{K}{E}$ with $f\mu=0$ such that,
when $G\in\catsa$ and $g\in\Mor{G}{E}$ satisfies $fg=0$, there exists a unique morphism $\overline{g}\in\Mor{G}{K}$
making the following diagram commutative
	 \begin{equation}\label{Diagram: kernel definition}
     \begin{tikzpicture}[auto, baseline=(current  bounding  box.center)]\matrix(m)[matrix of math nodes, column sep=4em, row sep=2em]
	 { K & E & F  \\
	 	{}  & G & {} \\};
 	\draw[->] (m-1-1) to node {\footnotesize{$\mu$}} (m-1-2); \draw[->] (m-1-2) to node {\footnotesize{$f$}} (m-1-3);
 		\draw[->, bend left] (m-1-1.north) to node {\footnotesize{$0$}} (m-1-3.north);
 	\draw[->] (m-2-2) to node {\footnotesize{$g$}} (m-1-2);
 	\draw[->] (m-2-2) to node [swap]{\footnotesize{$0$}} (m-1-3); \draw[->, dashed] (m-2-2) to node {\footnotesize{$\overline{g}$}} (m-1-1);
	 \end{tikzpicture}
	\end{equation}
	
A \textit{cokernel} of $f$ is a pair $(C,\pi)$, where $C\in\cats{A}$ and $\pi\in\Mor{F}{C}$ such that $\pi{}f=0$ and, whenever $G\in\catsa$ and $g\in\Mor{F}{G}$ satisfy $gf=0$, there exists a unique morphism $\overline{g}\in\Mor{C}{G}$ making the following diagram commutative
	 \begin{equation}\label{Diagram: cokernel definition}
     \begin{tikzpicture}[auto, baseline=(current  bounding  box.center)]\matrix(m)[matrix of math nodes, column sep=4em, row sep=2em]
	{ E & F & C \\
		{}  & G & {} \\};
	\draw[->] (m-1-1) to node {\footnotesize{$f$}} (m-1-2); \draw[->] (m-1-2) to node {\footnotesize{$\pi$}} (m-1-3);
	\draw[->, bend left] (m-1-1.north) to node {\footnotesize{$0$}} (m-1-3.north);
	\draw[->] (m-1-2) to node {\footnotesize{$g$}} (m-2-2);
		\draw[->] (m-1-1) to node [swap]{\footnotesize{$0$}} (m-2-2); \draw[->, dashed] (m-1-3) to node {\footnotesize{$\overline{g}$}} (m-2-2);
	\end{tikzpicture}
	\end{equation}
\end{definition}

\begin{remark}\label{Remark: Kernel and Cokernel properties are universal}
It is easy to see that kernel morphisms must be monomorphisms. We also note that the property of being a kernel of a morphism is universal.
That is, in Diagram~\ref{Diagram: kernel definition}, $(G,g)$ is a kernel for $f$ if and only if $\overline{g}$ is an isomorphism.
Similarly, cokernel morphisms are epimorphisms and cokernel objects are unique up to isomorphism.
	
It is not difficult to show that, in a category $\mathcal{A}$ where every morphism has a kernel and a cokernel, a morphism is a kernel
if and only if it is the kernel of its cokernel and is a cokernel if and only if it is the cokernel of its kernel.
Moreover, let the following diagram in $\mathcal{A}$ be commutative.
\begin{equation}\label{Diagram: Kernel-cokernel pairs are closed under isomorphisms}
\begin{tikzpicture}[auto, baseline=(current  bounding  box.center)]
\matrix(m)[matrix of math nodes, column sep=4em, row sep=2em]
{ \phantom{'}K & E\phantom{'} & C\phantom{'} \\
	\phantom{'}{F}  & E' & {G}\phantom{'} \\};
	\draw[->] (m-1-1) to node {\footnotesize{$\mu$}} (m-1-2); \draw[->] (m-1-2) to node {\footnotesize{$\pi$}} (m-1-3);
\draw[->] (m-1-1) to node { } (m-2-1); \draw[->] (m-1-2) to node { } (m-2-2);  \draw[->] (m-1-3) to node { } (m-2-3);
\draw[->] (m-2-1) to node {\footnotesize{$f$}} (m-2-2); \draw[->] (m-2-2) to node {\footnotesize{$g$}} (m-2-3);
\end{tikzpicture}
\end{equation}
Suppose $\mu$ is a kernel of $\pi$ and $\pi$ is a cokernel of $\mu$ and that all of the vertical arrows are isomorphisms.
Then, by the universal properties of kernels and cokernels, $f$ is a kernel of $g$ and $g$ is a cokernel of~$f$.
\end{remark}
\begin{definition}\label{Definition: Injectivity and Projectivity}
Let $\mathcal{A}$ be a category.
	
Suppose $\cats{M}$  is a class of monomorphisms in $\mathcal{A}$, closed under composition and such that every isomorphism in
$\catsa$ is in $\cats{M}$. An object $I\in\mathcal{A}$ is $\mathcal{M}$\textit{-injective} if, when given ${\mu}\in\classMofmorphisms{E}{F}$ and $f\in\Morunder{E}{I}{\mathcal{A}}$,  for objects $E,F\in\mathcal{A}$, there exists a morphism ${g}\in\Morunder{F}{I}{\mathcal{A}}$
making the following diagram commutative
\begin{equation}\label{Diagram: Injectivity}
\begin{tikzpicture}[auto, baseline=(current  bounding  box.center)]\matrix(m)[matrix of math nodes, column sep=3em, row sep=3em]
{ E & F  \\
	 I & {} \\};
\draw[->] (m-1-1) to node {\footnotesize{$\mu$}} (m-1-2);
\draw[->] (m-1-1) to node [swap]{\footnotesize{$f$}} (m-2-1);
\draw[->, dashed] (m-1-2) to node {\footnotesize{$g$}} (m-2-1);
\end{tikzpicture}
\end{equation}

Suppose $\cats{P}$ is a class of  epimorphisms in $\catsa$, closed under composition and such that every isomorphism in $\catsa$ is
in $\cats{P}$. An object $P\in\cats{A}$ is $\cats{P}$\textit{-projective} if, when given ${\pi}\in\classofmorphisms{P}{E}{F}$ and $f\in\Morunder{P}{F}{\mathcal{A}}$, for objects $E,F\in\mathcal{A}$, there exists a morphism  ${g}\in\Morunder{P}{E}{\mathcal{A}}$
making the following diagram commutative
\begin{equation}\label{Diagram: Projectivity}
\begin{tikzpicture}[auto, baseline=(current  bounding  box.center)]\matrix(m)[matrix of math nodes, column sep=3em, row sep=3em]
{ {} & {P} \\
E & F  \\};
\draw[->] (m-2-1) to node {\footnotesize{$\pi$}} (m-2-2);
\draw[->] (m-1-2) to node {\footnotesize{$f$}} (m-2-2);
\draw[->, dashed] (m-1-2) to node [swap]{\footnotesize{$g$}} (m-2-1);
\end{tikzpicture}
\end{equation}
We say $\mathcal{A}$ has \textit{enough $\mathcal{M}$-injectives} (resp., \textit{enough $\mathcal{P}$-projectives})
if, for every $E\in\mathcal{A}$, there exists an $\mathcal{M}$-injective object $I$ (resp., $\mathcal{P}$-projective object~$P$)
such that $\classMofmorphisms{E}{I}\neq\emptyset$ (resp., $\classofmorphisms{\mathcal{P}}{P}{E}\neq\emptyset$).
\end{definition}
\begin{remark}\label{Remark: Retractions of injectives and projectives}
Fix a category $\mathcal{A}$ and classes $\mathcal{M}$ and $\mathcal{P}$ of morphisms.
It is easy to see that any retract of an $\mathcal{M}$-injective object must be $\mathcal{M}$-injective and every retract
of a $\mathcal{P}$-projective object must be $\mathcal{P}$-projective. Moreover, if $I$ is $\mathcal{M}$-injective and
there exists a morphism $\mu\in\classMofmorphisms{I}{E}$, then $\mu$ must be a section of some retraction
$r\in\Morunder{E}{I}{\mathcal{A}}$. If $P$ is $\mathcal{P}$-projective and $\pi\in\classofmorphisms{P}{E}{P}$,
then $\pi$ must be a retraction.
\end{remark}

We will need the notion of kernels and cokernels when we talk about exact categories and the notion of injectives and
projectives when we discuss their global dimensions in Section~\ref{Section: Exact Categories and Global Dimension}.

\subsection{Operator Spaces}\label{Subsect:op-spaces}
In the background we will be using Ruan's Representation Theorem (\cite[Theorem 2.3.5]{ER1}).
Hence we will not distinguish between spaces arising as subspaces of $B(H)$, the bounded operators on a Hilbert space $H$,
and matrix normed spaces satisfying Ruan's axioms (see, for example, \cite[1.2.12]{BLM04}). However, when we refer to
\textit{operator spaces}, we will mean spaces of this type that we also assume to be complete. If $E$ is an operator space
we will write $x=[x_{i j}]\in M_n(E)$ to say that $x$ is an $n\times n$ matrix with entries $x_{i j}$, $i,j\in\{1,\ldots, n\}$ in $E$,
and $\norm{x}_n$ is the norm of $x$ in the Banach space $M_n(E)$, whose norm is inherited by the matrix norm of~$E$.

When $E,F$ are operator spaces, we will denote the operator space consisting of all completely bounded linear maps from
$E$ to ${F}$ by $\CB{E}{F}$. The completely bounded norm of an element $\phi\in\CB{E}{F}$ is~$\cbnorm{\phi}$.

By $A$ being an \textit{operator algebra} we will mean $A$ is a closed subalgebra of $B(H)$ of some Hilbert space $H$
such that $A$ contains a contractive approximate identity. By the Blecher--Ruan--Sinclair Theorem (\cite[Theorem 2.3.2]{BLM04}),
we need not distinguish between $A$ being an operator algebra and $A$ being an operator space which is also a Banach algebra
with a contractive approximate identity
such that the map $\hten{A}{A}\rightarrow{A}$ from the Haagerup tensor product of $A$ with itself to $A$ induced by the
multiplication on $A$ is completely contractive.

A complex algebra is said to be \textit{classically semisimple\/} if it is a direct sum of minimal right ideals and if it is finitely generated,
finitely many minimal right ideals suffice. Hence, combining the Artin--Wedderburn theorem with the Gelfand--Mazur theorem
it follows that a unital complex Banach algebra is classically semisimple if and only if it is
the direct sum of finitely many full matrix algebras over the complex numbers; in particular it is finite dimensional.

\begin{definition}\label{Definition: Operator module}
Let $A$ be an operator algebra. We say that a right $A$-module $E$ that is also an operator space is a \textit{right operator $A$-module}
if the map  $\hten{E}{A}\rightarrow{E}$, induced by the module action, is completely contractive.
\end{definition}

It follows from the associativity of the Haagerup tensor product that, for any operator space~$E$, $\hten EA$ is a right operator $A$-module;
this works analogous to \cite[3.1.5 (3)]{BLM04}.

\begin{definition}\label{Definition: Nondegenerate module}
Let $E$ be a right Banach $A$-module. If the linear span of elements of the form $x\cdot{a}$,
where $x\in E, a\in A$, is dense in $E$ then we say $E$ is \textit{non-degenerate}.
By Cohen's Factorisation Theorem \cite[Theorem~A.6.2]{BLM04}, this happens if and only if, for each $x\in E$,
there exist $x'\in E, a\in A$ such that $x=x'\cdot{a}$.
\end{definition}
We will always restrict ourselves to categories of non-degenerate operator modules, and unital modules if the algebra is unital.
If $E$ is an non-degenerate right operator $A$-module and $F$ is a closed subspace of $E$ then $F$ and $E/F$ become non-degenerate right operator $A$-modules when equipped with the standard induced operator space structures and module actions \cite[Section~3.1]{BLM04}.

For a Banach space $E$ and some $\gamma>0$ we will denote $\setst{x\in E}{\norm{x}<\gamma}$ by $E_{\norm{\cdot}<\gamma}$.
Recall that a linear map $\mapfromto{f}{E}{F}$ between Banach spaces is an open map if and only there exists $\gamma>0$ such that $f(E_{\norm{\cdot}<1})\supseteq F_{\norm{\cdot}<\gamma}$. If $E$ and $F$ are, moreover, operator spaces we say that $f$ is
\textit{completely open} if there exists some common $\gamma$ such that each amplification $f_n$ is open and for each $n\in\mathbb{N}$, $f_n\bigl(M_n(E)_{\norm{\cdot}<1}\bigr)\supseteq M_n(F)_{\norm{\cdot}<\gamma}$.
These maps can be characterised in the following way.

\begin{prop}\label{Proposition: Characterisation of CO maps}
Let $E$ and $F$ be operator spaces. Then $f\in\CB{E}{F}$ is completely open if and only if there exists $\lambda>1$
such that, for each $n\in\mathbb{N}$, every $y\in M_n(F)$ is equal to $f_n(x)$ for some $x\in M_n(E)$
with $\norm{x}_n\leq\lambda\norm{y}_n$.
\end{prop}	
For a completely open map $f$,
we will refer to $\lambda$ as in Proposition~\ref{Proposition: Characterisation of CO maps} as an \textit{openness constant} for~$f$.
\begin{example}\label{Example: Complete Quotient maps are open}
When $E$ is an operator space and $F$ is a closed subspace, then the canonical projection
$\pi\in\CB{E}{E/F}$ is a completely open map and any $\lambda>{1}$ is an openness constant.
\end{example}
Injectivity and self-duality of the Haagerup tensor product provide us with the following useful result.
\begin{lem}\label{lem:isom-embeddings}
Let $E$ and $F$ be operator spaces and let $u\in\hten EF$ be non-zero.
There exist bounded linear functionals $\alpha\in E^*$, $\beta\in F^*$ such that $(\alpha\otimes\beta)(u)\neq0$.
\end{lem}

\noindent
This follows immediately from the completely isometric embeddings
\[
\hten EF\hookrightarrow\hten{E^{**}}{F^{**}}\hookrightarrow(\hten{E^*}{F^*})^*
\]
given by $(x\otimes y)(\alpha\otimes\beta)=\alpha(x)\beta(y)$ for $x\in E$, $y\in F$, $\alpha\in E^*$ and $\beta\in F^*$,
see \cite[Proposition 9.2.5 and Theorem 9.4.7]{ER1}.

\section{The additive category of operator modules}\label{Section: Additive category of operator modules}
\noindent
It is well known how the notion of the global dimension for module categories extends to the setting of abelian categories with enough
injectives (or enough projectives). See, e.g., \cite{HiltStamm}, \cite{Mitch} or~\cite{O2000}.
In this section we will see to what extent the canonical additive category of non-degenerate (right) operator
$A$-modules over an operator algebra $A$ fails to be abelian. In order to remedy this by introducing an exact structure in the next section,
we need to study the kernels and cokernels in this category in detail.

We first recall the definition of an additive category.
\begin{definition}\label{Definition: Additive category}
A category $\mathcal{A}$ is \textit{additive} if $\mathcal{A}$ has a zero object; morphism sets have the structure of abelian groups;
composition is distributive with respect to this abelian group structure; and a product exists for each pair of objects.
\end{definition}
For $E,F\in\mathcal{A}$, where $\mathcal{A}$ is additive, we denote by $\Homunder{E}{F}{\mathcal{A}}$ the morphism set
equipped with the abelian group structure.
A functor $\mapfromto{\Ffunct}{\mathcal{A}}{\mathcal{B}}$ is \textit{additive} if, when $E,F\in\mathcal{A}$,
$\Ffunct(f+g)=\Ffunct(f)+\Ffunct(g)$ for all morphisms $f,g\in\Homunder{E}{F}{\mathcal{A}}$.
Note that, in an additive category, $E,F\in\mathcal{A}$ have a product $G\in\mathcal{A}$ if and only if $G$ is also their coproduct.
Moreover, this happens precisely when $G$ is a \textit{direct sum} of $E$ and $F$, that is, there exists a quintuplet
$(G,\iota_E,\iota_F,\pi_E,\pi_F)$, where $\iota_E\in\Hom{E}{G}, \iota_F\in\Hom{F}{G},$ and $\pi_E\in\Hom{G}{E}, \pi_F\in\Hom{G}{F}$, such that $\pi_E\iota_E=\idntyof{E},$ $\pi_F\iota_F=\idntyof{F}$ and $\iota_E\pi_E+\iota_F\pi_F=\idntyof{G}$.

For the remainder of this section we fix an operator algebra~$A$.
We will use $\OMODAinf$ to denote the category whose objects are the non-degenerate right operator $A$-modules and whose morphisms
are the completely bounded $A$-module maps. Instead of $\OMODCCinf$ we write~$\Opinf$.
For $E,F\in\OMODAinf$, we denote the set $\Morunder{E}{F}{\OMODAinf}$ by $\CBmod{A}{E}{F}$, and by $\CB EF$ if $A=\CC$.
It is clear that this is an abelian group.

\begin{prop}\label{Proposition: The Category OMODAinf is additive}
Let $A$ be an operator algebra. The category $\OMODAinf$ is additive.
\end{prop}
The only part of Proposition~\ref{Proposition: The Category OMODAinf is additive} that is perhaps not immediately apparent is the existence
of a product for each pair of objects in $\OMODAinf$.  For any operator spaces $E_1, E_2,$ let $E_1\oplus E_2$ be the vector space direct sum
equipped with the norm $\norm{(x,y)}=\norm{x}+\norm{y}$ for each $x\in E_1, y\in E_2$. For each $n\in\mathbb{N}$, the obvious identifications
$M_n(E_1\oplus E_2)\cong M_n(E_1)\oplus M_n(E_2)$ yield an operator space structure on $E_1\oplus E_2$.
Moreover, if $E_1, E_2\in\OMODAinf,$ then we equip $E_1\oplus E_2$ with the module action
$(x,y)\cdot{a}:=(x\cdot{a},y\cdot{a})$ for all $x\in E_1, y\in E_2, a\in A$.
With this, we have $E_1\oplus E_2\in\OMODAinf$ and the quintuplet $(E_1\oplus E_2,\iota_{1},\iota_{2},\pi_{1},\pi_{2})$,
where for each $i\in\{1,2\}$, $\mapfromto{\iota_{i}}{E_i}{E_1\oplus E_2}$ denotes the inclusion and
$\mapfromto{\pi_{i}}{E_1\oplus E_2}{E_i}$ is the projection, makes $E_1\oplus E_2$ a direct sum of $E_1$ and~$E_2$.

The next result is elementary.
\begin{prop}\label{Proposition: Monomorphisms and Epimorphisms in OMODAinf}
Let $E,F\in\OMODAinf$. For every $f\in\CBmod{A}{E}{F}$, we have
\begin{enumerate}[label=\upshape(\roman*), nosep]
		\item $f$ is a monomorphism in $\OMODAinf$ if and only if it is injective.
		\item $f$ is an epimorphism in $\OMODAinf$ if and only if it has dense range.
\end{enumerate}
\end{prop}

Kernels and cokernels are fundamental tools in our approach.
\begin{prop}\label{Proposition: Every morphism in OMODAinf has a kernel and a cokernel}
Let $A$ be an operator algebra. Suppose $E,F\in\OMODAinf$ and $f\in\CBmod{A}{E}{F}$. Then $f$ has a kernel and a cokernel.
\end{prop}
\begin{proof}
We first deal with the kernel of~$f$. Let $K$ be $\ker{f}=f^{-1}(0)$; this is clearly an object in $\OMODAinf$.
Let $\mapfromto{\mu}{K}{E}$ be the inclusion. Then $\mu\in\CBmod AEF$ and $f\mu=0$.
Suppose $G\in\OMODAinf$ and there exists $g\in\CBmod{A}{G}{E}$ with $fg=0$.
Then $g(G)\subseteq{K},$ so we let $\overline{g}\in\CBmod{A}{G}{K}$ just be the $A$-module map~$g$.
As morphisms we have $g=\mu\overline{g}$. That this is the only such morphism making Diagram~\eqref{Diagram: kernel definition} commutative,
follows from the fact that $\mu$ is a monomorphism.
		
To prove that $f$ has a cokernel, we let $C=F/{\overline{f(E)}}$ and $\pi\in\CBmod{A}{F}{C}$ be the canonical projection.
Obviously, $\pi{}f=0$. Suppose there  exist $G\in\OMODAinf$ and $g\in\CBmod{A}{F}{G}$ such that $gf=0$.
For each $y\in F,$ let $\overline{g}(\pi(y))=g(y)$. If $\pi(y)=0$ then $y\in\overline{F(E)}$ and, by continuity, $g(y)=0$.
Hence $\mapfromto{\overline{g}}{C}{G}$ is a well-defined map and is clearly an $A$-module map.
In fact, it is completely open with openness constant~$\lambda>1$ since $M_n(C)\cong M_n(F)/M_n(\overline{F(E)})$
\cite[1.2.14]{BLM04}.
For each $n\in\mathbb{N}$ and $c\in M_n(C)$ choose $y\in M_n(F)$ such that $\pi_{n}(y)=c$ and $\norm{y}_n\leq\lambda\norm{c}_n$.
Then $\norm{\overline{g}(c)}_n=\norm{g(y)}_n\leq\cbnorm{g}\lambda\norm{c}_n$.
Hence $\overline{g}$ is completely bounded, i.e., $\overline{g}\in\CBmod{A}{C}{G}$ such that $\overline{g}\pi=g$.
Lastly, we note since $\pi$ is an epimorphism, $\overline{g}$ is the only morphism $C\rightarrow{G}$ making
Diagram~\eqref{Diagram: cokernel definition} commutative.
\end{proof}
\begin{prop}\label{Proposition: Description of kernels in OMODAinf}
Let $A$ be an operator algebra. Let $E,F\in\OMODAinf$ and $\mu\in\CBmod{A}{E}{F}$.
Then $\mu$ is the kernel of some morphism $\pi\in\CBmod{A}{F}{G}$,
$G\in\OMODAinf$ if and only if $\mu$ has closed range and is an isomorphism in $\OMODAinf$
(that is, a completely bounded bijective $A$-module map, with completely bounded inverse) when viewed as a map onto its range.
\end{prop}
\begin{proof}
By Remark~\ref{Remark: Kernel and Cokernel properties are universal} and Proposition~\ref{Proposition: Every morphism in OMODAinf has a kernel and a cokernel} we need only consider the case where $G=F/\overline{\mu(E)}$ and $\pi$ is the canonical quotient mapping.
Let $K=\ker{\pi}$ and $\iota\in\CBmod{A}{K}{F}$ be the inclusion map. Then $\mu$ is a kernel of $\pi$ if and only if there exists an isomorphism
$\phi\in\CBmod AKE$ making the following diagram commutative.
		\begin{equation}\label{Diagram: Description of kernels in OMODAinf}
        \begin{tikzpicture}[auto, baseline=(current  bounding  box.center)]\matrix(m)[matrix of math nodes, column sep=4em, row sep=2em]
		{ E & F & G  \\
			{}  & K & {} \\};
		\draw[->] (m-1-1) to node {\footnotesize{$\mu$}} (m-1-2); \draw[->] (m-1-2) to node {\footnotesize{$\pi$}} (m-1-3);
		\draw[->, bend left] (m-1-1.north) to node {\footnotesize{$0$}} (m-1-3.north);
		\draw[->] (m-2-2) to node {\footnotesize{$\iota$}} (m-1-2);
		\draw[->] (m-2-2) to node [swap]{\footnotesize{$0$}} (m-1-3); \draw[->, dashed] (m-2-2) to node {\footnotesize{$\phi$}} (m-1-1);
		\end{tikzpicture}
		\end{equation}
Note that $\overline{\mu(E)}=\iota(K)$. Suppose $\mu$ has closed range and is an isomorphism in $\OMODAinf$
when considered as a map onto its range. Then $\mu(E)=\iota(K)$ and we simply let $\phi$ be the completely bounded inverse
$\mu(E)\rightarrow{E}$ composed with~$\iota$. So $\mu$ is a kernel.
		
On the other hand, if $\mu$ is the kernel of $\pi$ then $\phi$ exists. Since $\mu=\iota\circ\phi^{-1}$ and $\iota$ is an isometry,
$\mu(E)$ is closed and we are done.
\end{proof}
\begin{prop}\label{Proposition: Description of cokernels in OMODAinf}
Let $A$ be an operator algebra. Let $E,F\in\OMODAinf$ and $\pi\in\CBmod{A}{E}{F}$.
Then $\pi$ is the cokernel of some morphism $\mu\in\CBmod{A}{G}{E}$, $G\in\OMODAinf$ if and only if $\pi$ is completely open.
\end{prop}
\begin{proof}
By Remark~\ref{Remark: Kernel and Cokernel properties are universal} and
Proposition~\ref{Proposition: Every morphism in OMODAinf has a kernel and a cokernel}
we need only look at the case where $G=\ker{\pi}$ and $\mu$ is the inclusion mapping.
Let $C=E/{\mu(G)}$ and $g\in\CBmod{A}{E}{C}$ be the quotient map.
Then $\pi$ is cokernel for $\mu$ if and only if there exists an isomorphism $\mapfromto{\phi}{F}{C}$ making the following diagram commutative.
		\begin{equation}\label{Diagram: Description of cokernels in OMODAinf}
        \begin{tikzpicture}[auto, baseline=(current  bounding  box.center)]\matrix(m)[matrix of math nodes, column sep=4em, row sep=2em]
		{ G & E & F \\
			{}  & C & {} \\};
		\draw[->] (m-1-1) to node {\footnotesize{$\mu$}} (m-1-2); \draw[->] (m-1-2) to node {\footnotesize{$\pi$}} (m-1-3);
		\draw[->, bend left] (m-1-1.north) to node {\footnotesize{$0$}} (m-1-3.north);
		\draw[->] (m-1-2) to node {\footnotesize{$g$}} (m-2-2);
		\draw[->] (m-1-1) to node [swap]{\footnotesize{$0$}} (m-2-2); \draw[->, dashed] (m-1-3) to node {\footnotesize{$\phi$}} (m-2-2);
		\end{tikzpicture}
		\end{equation}
Suppose that $\pi$ is completely open and hence surjective. Note that, if $x\in E$ is such that $\pi(x)=0$,
then $x\in\mu(G)$ and $g(\mu(x))=0$, so the map $\mapfromto{\phi}{F}{C}$, $\pi(x)\mapsto g(x)$ is well defined.
As $\pi, g$ are $A$-module maps so is~$\phi$. For any $n\in\mathbb{N}$ and $y\in M_n(F)$, we have by
Proposition~\ref{Proposition: Characterisation of CO maps} that
$\norm{\phi_n(y)}_n=\norm{g_n(x)}_n\leq\lambda\cbnorm{g}\norm{y}_n$ for some $x\in E$ and openness constant $\lambda>1$.
So $\phi\in\CBmod{A}{F}{C}$. A similar argument (using the fact that $g$ is completely open) gives that  there is a morphism $\psi\in\CBmod{A}{C}{F}$ defined by setting $\psi(g(x))=\pi(x)$ for any $x\in E$.
By definition, $\phi\pi=g$ and it is clear that $\psi$ is the inverse of $\phi$, thus $\phi$ is an isomorphism.
It follows that $\pi$ is a cokernel of~$\mu$.
	
Conversely, suppose there exists an isomorphism $\phi$ making Diagram~\eqref{Diagram: Description of cokernels in OMODAinf} commutative.
Let $n\in\mathbb{N}$ and $y\in M_n(F)$. By Proposition~\ref{Proposition: Characterisation of CO maps}, there exists $x\in M_n(E)$ such that $g_n(x)=\phi_n(y)$ and $\norm{x}_n\leq\lambda\norm{\phi_n(y)}_n\leq\lambda\cbnorm{\phi}\norm{y}_n$, where $\lambda$ is an
openness constant for~$g$.
Moreover, the commutativity of Diagram~\eqref{Diagram: Description of cokernels in OMODAinf} gives that
$\pi_n(x)=\phi_n^{-1}g_n(x)=\phi_n^{-1}\phi^{}_n(y)=y$ and by Proposition~\ref{Proposition: Characterisation of CO maps}, $\pi$ is completely open.
\end{proof}
An additive category $\mathcal{A}$ is called \textit{abelian} if
every morphism in $\mathcal{A}$ has both a kernel and a cokernel;
every monomorphism in $\mathcal{A}$ is the kernel of its cokernel;
and every epimorphism in $\mathcal{A}$ is the cokernel of its kernel.
Comparing the descriptions above of the kernels and cokernels in $\OMODAinf$ with the monomorphisms and epimorphisms in that category, it is perhaps obvious that $\OMODAinf$ fails to be abelian. Here is the probably simplest example for $A=\CC$.
Let $E$ be an infinite-dimensional Banach space. Recall that the identity map on $E$ is completely contractive when considered as a map $\mapfromto{f}{\textup{Max}(E)}{\textup{Min}(E)}$, where ${\textup{Max}(E)}$ is $E$ equipped with its maximal quantization and
${\textup{Min}(E)}$ is its minimal quantization. Then $f\in\CBmod{\mathbb{C}}{\textup{Max}(E)}{\textup{Min}(E)}$ is surjective
and hence an epimorphism in $\OMODinf{\mathbb{C}}$.
However it cannot be completely open, that is, an isomorphism as $E$ is infinite dimensional (see, e.g., \cite[Corollary 3.9]{P03}).

As a result, the homological algebra which is highly developed in abelian categories is not directly available in $\OMODAinf$.
In the next section we show how we can get around this issue  by introducing exact categories. The following two results will be essential.
\begin{prop}\label{Proposition: Pullbacks and cokernels in OMODAinf}
Let $E,F,G\in\OMODAinf$ with $f_E\in\CBmod{A}{E}{G}$ and $f_F\in\CBmod{A}{F}{G}$.
	\begin{enumerate}[label=\upshape(\roman*), leftmargin=*]
	\item There exist $L\in\OMODAinf$ and $\ell_E\in\CBmod{A}{L}{E}$, $\ell_F\in\CBmod{A}{L}{F}$ making the following diagram a pullback square.
	\begin{equation}\label{Diagram: Pullback in OMODAinf}
	\begin{tikzpicture}[auto, baseline=(current  bounding  box.center)]
	\matrix(m)[matrix of math nodes, column sep=4em, row sep=3em]{ L & F \\ E & G\\	};
	\draw[->,dashed] (m-1-1) to node {\footnotesize{$\ell_F$}} (m-1-2); \draw[->,dashed] (m-1-1) to node [swap]{\footnotesize{$\ell_E$}} (m-2-1);
	\draw[->] (m-1-2) to node {\footnotesize{$f_F$}} (m-2-2); \draw[->] (m-2-1) to node [swap]{\footnotesize{$f_E$}} (m-2-2);
	\end{tikzpicture}
	\end{equation}
\item	For any such pullback square, if $f_E$ is a cokernel map, then so is $\ell_F$.
	\end{enumerate}
\end{prop}
\begin{proof}
	(i)	Let $L=\setst{(x,y)\in E\oplus F}{f_E(x)=f_F(y)}$.
Then $L$ is a closed submodule of $E\oplus F$ so inherits the operator $A$-module structure of $E\oplus F$.
Let $\ell_F$ and $\ell_E$ be the restrictions to $L$ of the canonical projections $\mapfromto{\pi_F}{E\oplus{F}}{F}$
and $\mapfromto{\pi_E}{E\oplus{F}}{E}$, respectively. By definition of~$L$,
Diagram~\eqref{Diagram: Pullback in OMODAinf} is commutative.
		
If there exist $L'\in\OMODAinf$ and $\ell'_E\in\CBmod{A}{L'}{E},$  $\ell'_F\in\CBmod{A}{L'}{F}$
such that $f_F\ell'_F=f_E\ell'_E$, then, by the universal property of products, there exists a unique $\phi\in\CBmod{A}{L'}{L}$ such that $\ell'_E=\pi_{E}\phi$  and $\ell'_F=\pi_{F}\phi$ and it is clear that $\phi(L')\subseteq{L}$; hence $L$ must make
Diagram~\eqref{Diagram: Pullback in OMODAinf} a pullback square.

	(ii)  By Remark~\ref{Remark: Kernel and Cokernel properties are universal} and the universal property of pullbacks,
the result holds if and only if it holds for the pullback square defined in~(i).
Suppose that $f_E$ is a cokernel map. Proposition~\ref{Proposition: Description of cokernels in OMODAinf}
tells us that $f_E$ is completely open and we are done if $\ell_F$ is completely open.
	
Let $\lambda$ be an openness constant for $f_{E}$ and set $\lambda'=\max{\{\lambda\cbnorm{f_F}, 1\}}$.
For $n\in\mathbb{N}$ and $y\in M_n(F)$, we have $(f_F)_n(y)\in M_n(G)$ and,
by Proposition~\ref{Proposition: Characterisation of CO maps}, there exists $x\in M_n(E)$ such that
$(f_E)_n(x)=(f_F)_n(y)$ (hence $(x, y)\in L$) with $\norm{x}_n\leq\lambda\norm{(f_F)_n(y)}_n\leq\lambda'\norm{y}_n$.
We note that $(\ell_F)_n(x,y)=y$ with $\norm{(x,y)}_n\leq2\lambda'\norm{y}_n$.
Proposition~\ref{Proposition: Characterisation of CO maps} tells us $\ell_F$ is completely open.
\end{proof}
\begin{prop}\label{Proposition: Pushouts and kernels in OMODAinf}
Let $E,F,G\in\OMODAinf$ with $f_E\in\CBmod{A}{G}{E}$ and $f_F\in\CBmod{A}{G}{F}$.
	\begin{enumerate}[label=\upshape(\roman*), leftmargin=*]
    \item There exist $C\in\OMODAinf$ and $h_E\in\CBmod{A}{E}{C}, h_F\in\CBmod{A}{F}{C}$ making the following diagram a pushout square.
	\begin{equation}\label{Diagram: Pushout in OMODAinf}
	\begin{tikzpicture}[auto, baseline=(current  bounding  box.center)]
	\matrix(m)[matrix of math nodes, column sep=4em, row sep=3em]{ G & E \\ F & C\\	};
	\draw[->] (m-1-1) to node {\footnotesize{$f_E$}} (m-1-2); \draw[->] (m-1-1) to node [swap]{\footnotesize{$f_F$}} (m-2-1);
	\draw[->,dashed] (m-1-2) to node {\footnotesize{$h_E$}} (m-2-2); \draw[->,dashed] (m-2-1) to node [swap]{\footnotesize{$h_F$}} (m-2-2);
	\end{tikzpicture}
	\end{equation}
    \item	For any such pushout square, if $f_E$ is a kernel map, then so is $h_F$.
\end{enumerate}
\end{prop}
\begin{proof}
		(i) Let $H=\setst{(f_E(z),-f_F(z))}{z\in G}$ and $C= {E\oplus{F}}/\overline{H}$.
Let $h_E$ be the composition of the embedding $\mapfromto{\iota_E}{E}{E\oplus{F}}$ with the canonical projection $\mapfromto{\pi}{E\oplus{F}}{C}$ and $h_F=\pi\iota_F$, where $\iota_F$ is the embedding ${F}\rightarrow{E\oplus{F}}$.
Clearly $h_E\in\CBmod{A}{E}{C}$ and $h_F\in\CBmod{A}{F}{C}.$ For any $z\in G,$ $(f_E(z),0) - (0, f_F(z))=(f_E(z), -f_F(z))\in H$; this means that Diagram~\eqref{Diagram: Pushout in OMODAinf} is commutative.
		
Suppose there exists $C'\in\OMODAinf$ and $h'_E\in\CBmod{A}{E}{C'}$ and $h'_F\in\CBmod{A}{F}{C'}$ with $ h'_Ef_E=h'_Ff_F$.
By the universal property of coproducts, there exists $\phi\in\CBmod{A}{E\oplus{F}}{C'}$ such that $\phi\iota_E=h'_E$ and $\phi\iota_F=h'_F$.
For Diagram~\eqref{Diagram: Pushout in OMODAinf} to be a pushout square it remains to show that there exists $\tau\in\CBmod{A}{C}{C'}$ such that $\tau{}h_E=h_E'$ and $\tau{}h_F=h_F'$.
		Suppose $(x,y)\in H$, then there exists $z\in G$ such that $x=f_E(z)$ and $y=-f_F(z)$.
Therefore,
\[
\phi(x,y)=h_E'(f_E(z))+h_F'(-f_F(z))=h_E'(f_E(z))-h_E'(f_E(z))=0.
\]
That is, the $A$-module map $\mapfromto{\tau}{C}{C'},\, \pi(x,y)\mapsto\phi(x,y)$ is well defined.
Let $c\in M_n(C)$ and $\lambda>1$. As in Example~\ref{Example: Complete Quotient maps are open},
$\pi$ is completely open and $\lambda$ is an openness constant for $\pi$.
Therefore, there exists $(x,y)\in M_n(E\oplus{F})$ such that $c=\pi_n(x,y)$ with $\norm{\tau{(c)}}_n=\norm{\phi(x,y)}_n\leq\cbnorm{\phi}\norm{(x,y)}_n\leq\lambda \cbnorm{\phi}\norm{c}_n$.
Hence $\tau$ is completely bounded.
		
(ii)  By Remark~\ref{Remark: Kernel and Cokernel properties are universal} and the universal property of pushouts,
the result holds if and only if it holds for the pushout square defined in~(i). Suppose that $f_E$ is a kernel map in $\OMODAinf$.
That is, $f_E(G)$ is closed in $E$ and there exists $g\in\CBmod{A}{f_E(G)}{G}$ such that $gf_E=\idntyof{G}$ and $f_Eg=\idntyof{f_E(G)}$.
We will show that $h_F$ is a kernel map too.
		
Suppose we have a sequence $(f_E(z_n), -f_F(z_n))_{n\in\mathbb{N}}$ in $H$ with limit $(x,y)\in E\oplus F$.
By continuity, $g(x)$ is the limit of $(z_n)_{n\in\mathbb{N}}=(gf_E(z_n))_{n\in\mathbb{N}}$ and $y=-f_F(g(x))$.
That is, $(x,y)=(f_E(g(x)), -f_F(g(x)))\in H$. Therefore $H$ is closed and $C=E\oplus{F}/H$.
		
		For $h_F$ to be a kernel map, we need that $h_F(F)$ is closed in~$C$.
Let $(h_F(y_n))_{n\in\mathbb{N}}\subseteq C$ and $c\in C$ be such that $\norm{h_F(y_n)-c}\to{0}$.
There exist $x\in E$, $y\in F$ with $c=\pi(x,y)$ is the limit of $(\pi(0,y_n))_{n\in\mathbb{N}}$;
that is, $\norm{\pi(-x, y_{n}-y)-\pi(0,y_n)}\to0$.
Because $H$ is closed in $E\oplus{F}$ there must exist a sequence $(z_n)\in G$ such that $(f_E(z_n),-f_F(z_n))\in H$ with $\norm{(-x-f_E(z_n), y_n-y+f_F(z_n))}\to0$.
In particular, $\norm{-x-f_E(z_n)}\to0$ and, as $f_E(G)$ is closed in $E$, there exists some $z\in G$ with $f_E(z)=-x$.
By continuity, $z_n=g(f_E(z_n))\to z$ and $f_F(z_n)\to f_F(z)$.
Set $y'=y-f_F(z)$, then $\norm{y_n-y'}\to0$. Therefore
\[
h_F(y')=\pi(0,y')=\pi((x,y)+(f_E(z), -f_F(z)))=\pi(x,y)
\]
and hence, $\pi(x,y)\in h_F(F)$.
		
Note that $h_F$ is injective. Indeed, if $h_F(y)=0,$ then there exists $z\in G$ such that $(0,y)=(f_E(z),-f_F(z))$
but $f_E$ is injective so $z=0$ and therefore $y=0$. So we certainly have an $A$-module map
$\mapfromto{\ell}{h_F(F)}{F}$ defined by $\ell(h_F(y))=y$.
We are done if $\ell$ is completely bounded.
		
Note that, for each $z\in M_n(G),$ we have
\begin{equation}\label{Equation: Equation 1 in kernel-pushout-OMODAinf proof}
		\norm{(f_F)_n(z)}_n\leq\cbnorm{f_F}\norm{z}_n\leq\cbnorm{f_F}\,\cbnorm{g}\,\norm{(f_E)_n(z)}_n.
\end{equation}
If $f_F=0$, the result is obvious, so we can suppose otherwise. Then
Equation~\eqref{Equation: Equation 1 in kernel-pushout-OMODAinf proof} tells us that $\norm{f_E(z)}_n\geq K{\norm{f_F(z)}_n}$ where $K=\min{\{(\cbnorm{f_F}\,\cbnorm{g})^{-1}, 1\}}$.
Recall that for each $y\in M_n(F)$, $\norm{(h_F)_n(y)}_n=\inf\norm{(0,y)+((f_E)_n(z), -(f_F)_n(z))}_n$,
where the infimum is over all $z\in M_n(G)$.
Then, for each $n\in\mathbb{N}$, there exists $z\in M_n(G)$ such that
\[
\begin{split}
\norm{(h_F)_n(y)}_n &\geq\frac{1}{2}\norm{(0,y)+((f_E)_n(z),-(f_F)_n(z))}_n=\frac{1}{2}\norm{((f_E)_n(z),y-(f_F)_n(z))}_n \\
		&=\frac{1}{2}\bigl(\norm{(f_E)_n(z)}_n+\norm{y-(f_F)_n(z)}_n\bigr) \\
        &\geq\frac{K}{2}\bigl(\norm{(f_F)_n(z)}_n+\norm{y-(f_F)_n(z)}_n\bigr)\geq\frac{K}{2}\norm{y}_n. 	
\end{split}
\]
Therefore $\norm{\ell_n((h_F)_n(y))}_n=\norm{y}_n\leq\frac{2}{K}\norm{(h_F)_n(y)}_n$ for all~$n$ and $\ell$ is completely bounded.
\end{proof}
We are now ready to embark on setting up the new framework for homological algebra.

\section{Exact Categories and Global Dimension}\label{Section: Exact Categories and Global Dimension}
\noindent
Exact categories are additive categories equipped with a class of kernel-cokernel pairs that satisfy certain axioms
(see Definition~\ref{Definition: Exact Category} below). In this section we show how exact categories provide a framework
to arrive at a more general version of the notion of a global dimension for an abelian category.
Our main source for exact categories is the survey article of B\"{u}hler~\cite{B08}.
We show that for a general operator algebra $A$, $\OMODAinf$ has a canonical exact structure.
\begin{definition}
Suppose we have a pair of composable morphisms
	\[\begin{tikzpicture}[auto] \node (A) {$K$};
	\node (B) [right=1.1cm of A] {$E$};
	\node (C) [right=1.1cm of B] {$C$};
	
	\draw[->] (A) to node {\footnotesize{$\mu_K$}} (B);
	\draw[->] (B) to node {\footnotesize{$\pi_C$}} (C);
	\end{tikzpicture}\]
in an additive category~$\mathcal{A}$, where $\mu_K$ is a kernel of $\pi_C$ and $\pi_C$ is a cokernel of~$\mu_K$.
Then we say that $(\mu_K, \pi_C)$ is a \textit{kernel-cokernel pair}.
Suppose $\mathcal{E}$ is a fixed class of kernel-cokernel pairs in~$\catsa$.
Then a morphism $\mu$ is called an \textit{admissible monomorphism}\index{Admissible! monomorphism} if there exists a morphism $\pi$ such that $(\mu,\pi)\in\mathcal{E}$. A morphism $\pi$ is called an \textit{admissible epimorphism}\index{Admissible! epiomorphism} if there exists a morphism $\mu$ such that $(\mu,\pi)\in\mathcal{E}$.
In diagrams admissible monomorphisms (resp., epimorphisms) will be depicted by $\rightarrowtail$  (resp., $\twoheadrightarrow$).
\end{definition}
We define exact categories using the axioms of an exact structure given in~\cite{B08}.
In Section~2 of that paper, B\"{u}hler shows these axioms are equivalent to Quillen's axioms (\cite{Quill73}).
\begin{definition}\label{Definition: Exact Category}
An \textit{exact structure} on $\mathcal{A}$ is a class $\Ex$ of kernel-cokernel pairs
which is closed under isomorphisms and satisfies the following axioms:
\begin{enumerate}
\item[[E$0$\!\!\!]] \ $\forall\ E\in\catsa\colon \idntyof{E}$ is an admissible monomorphism;
\item[[E$0^{\textup{op}}$\!\!\!]] \ $\forall\ E\in\catsa\colon \idntyof{E}$ is an admissible epimorphism;
\item[[E$1$\!\!\!]] \ the class of admissible monomorphisms is closed under composition;
\item[[E$1^{\textup{op}}$\!\!\!]] \ the class of admissible epimorphisms is closed under composition;
\item[[E$2$\!\!\!]] \ the pushout of an admissible monomorphism along an arbitrary morphism exists and yields
                    an admissible monomorphism;
\item[[E$2^{\textup{op}}$\!\!\!]] \ the pullback of an admissible epimorphism along an arbitrary morphism exists and yields
                                    an admissible epimorphism.
\end{enumerate}
An \textit{exact category} is a pair ($\catsa, \Ex$) consisting of an additive category~$\catsa$ and an exact structure $\Ex$ on~$\catsa$.
\end{definition}
\begin{definition}\label{Definition: split kernel-cokernel pair}
For an additive category $\mathcal{A}$, a kernel-cokernel pair
\[
\begin{tikzpicture}[auto] \node (A) {$K$};
	\node (B) [right=1.1cm of A] {$E$};
	\node (C) [right=1.1cm of B] {$C$};
	
	\draw[->] (A) to node {\footnotesize{$\mu_K$}} (B);
	\draw[->] (B) to node {\footnotesize{$\pi_C$}} (C);
\end{tikzpicture}
\]
is \textit{split} if there exist morphisms $\mu_C\in\Hom{C}{E}$ and $\pi_K\in\Hom{E}{K}$ such that the quintuplet
$(E,\mu_K,\mu_C,\pi_K,\pi_C)$ makes $E$ a direct sum.
\end{definition}
\begin{example}\label{Example:Min class of kernel-cokernel pairs}
The class of all split kernel-cokernel pairs forms an exact structure on any additive category.
We denote this class $\Exmin$ and note that $\Exmin$ will be a substructure of any other exact structure placed on~$\mathcal{A}$.
This is trivial to show, but has useful consequences
(see Section~\ref{Section: Relative Homological Algebra for Operator Modules over Operator Algebras} below).
\end{example}
\begin{example}\label{Example: In abelian categories the kernel-cokernel pairs form exact structure}
The class of all kernel-cokernel pairs in an additive category $\mathcal{A}$ is denoted~$\Exmax$.
If $\mathcal{A}$ is an abelian category then
$(\mathcal{A},\Exmax)$ is an exact category and $\Exmax$ is the class of all short exact sequences in~$\mathcal{A}$.
\end{example}

For a general additive category $\mathcal{A}$, the class $\Exmax$ is closed under isomorphisms
(Remark~\ref{Remark: Kernel and Cokernel properties are universal}) and satisfies [E$0$] and  [E$0^{\textup{op}}$].
We introduce some conditions that ensure that this class forms an exact structure on~$\mathcal{A}$.

\begin{definition}\label{Definition: Quasi-abelian category}
Let $\mathcal{A}$ be an additive category. We say $\mathcal{A}$ is a \textit{quasi-abelian category} if:
\begin{enumerate}[nosep,label=\upshape(\roman*)]
\item each morphism in $\mathcal{A}$ has a kernel and a cokernel;
\item the class $\Exmax$ of all kernel-cokernel pairs satisfies conditions [E2] and [E$2^\textup{op}$]
from Definition~\ref{Definition: Exact Category}.
\end{enumerate}
\end{definition}

The following can be proven by diagram chasing, from the definitions of kernels, cokernels, pushouts and pullbacks.
(See \cite[Proposition~4.4]{B08}.)
\begin{prop}\label{Proposition: Quasi-abelian categories have max exact structure}\textup{(\cite[Proposition 1.1.7]{Schneiders1999})}
Let $\mathcal{A}$ be a quasi-abelian category. Then $\Exmax$, the class of all kernel-cokernel pairs, satisfies conditions
\textup{[E1]} and \textup{[E$1^\textup{op}$]} from Definition~\ref{Definition: Exact Category} and thus forms an exact structure.
\end{prop}
\begin{theorem}\label{Theorem: Emax in OMODAinf forms an exact structure}
Let $A$ be an operator algebra. The class of all kernel-cokernel pairs in $\OMODAinf$ forms an exact structure on $\OMODAinf$.
\end{theorem}
\begin{proof}
By Propositions~\ref{Proposition: Every morphism in OMODAinf has a kernel and a cokernel},
\ref{Proposition: Pullbacks and cokernels in OMODAinf}
and~\ref{Proposition: Pushouts and kernels in OMODAinf}, $\OMODAinf$ is quasi-abelian.
The result then follows from Proposition~\ref{Proposition: Quasi-abelian categories have max exact structure}.
Alternatively, one can check \textup{[E1]} and \textup{[E$1^\textup{op}$]} `by hand' \cite[Section~3.4]{Rosb2020}.
\end{proof}
In an abelian category, the short exact sequences are the smallest building blocks for homological algebra
which then emerges from long exact sequences and so-called diagram lemmas. The kernel-cokernel pairs take
the role of short exact sequences, and the axioms of an exact category entail enough of the diagram lemmas to
build a viable theory. This is demonstrated in detail in~\cite{B08} and successfully applied to a functional analytic setting
in~\cite{B11}.

Another fundamental tool is the concept of an exact functor.
\begin{definition}\label{Definition: Exact functor}
An additive functor $\mapfromto{\Ffunct}{(\mathcal{A},\Ex_{1})}{(\mathcal{B}, \Ex_{2})}$ between two exact categories is \textit{exact}
if $\Ffunct(\mathcal{\Ex_{1}})\subseteq{\Ex_{2}}$.
\end{definition}

We will note for later (see Section~\ref{Section: Relative Homological Algebra for Operator Modules over Operator Algebras})
how exact functors are useful for inducing alternate exact structures on an exact category.

\begin{prop}\label{Proposition: Using exact functors to get exact structures}
Let $\mapfromto{\Ffunct}{(\mathcal{A},\Ex_{1})}{(\mathcal{B}, \Ex_{2})}$ be an exact functor between exact categories.
Suppose there is another exact structure $\Ex'_{2}$ on~$\mathcal{B}$. Then
\[
\Ex'_{1}=\setst{(\mu,\pi)\in\Ex_{1}}{\mathcal{(\Ffunct\mu,\Ffunct\pi)\in\Ex'_{2}}}
\]
forms an exact structure on~$\mathcal{A}$.
\end{prop}

\noindent
Proposition~\ref{Proposition: Using exact functors to get exact structures} is \cite[Proposition~7.3]{Heller1958};
the proof also follows easily from \cite[Proposition~5.2]{B08}.
We will call the exact structure $\Ex'_{1}$ appearing in this way a \textit{relative exact structure\/}
as it involves constraints relative to another exact structure.

Injectivity and projectivity can be described using the notion of exact functors.
For any object $E$ in an additive category $\mathcal{A}$ we have the following contravariant functor
\[
\begin{split}
\Hom{-}{E}\colon{\mathcal{A}}&\xrightarrow{\hphantom{fgl}}\Abcat\\
	F&\longmapsto\Homunder{F}{E}{\mathcal{A}}\\
	\Homunder{F}{G}{\mathcal{A}}\ni{f}&\longmapsto f^*=\Hom{f}{E}
\end{split}
\]
where $f^*(g)=gf$  for $g\in\Homunder{G}{E}{\mathcal{A}}$ and $\Abcat$ is the category of abelian groups.
There is also the similarly defined covariant functor $\Hom{E}{-}$, where $f_*=\Hom{E}{f}$ is given by post-composition
of~$f$. We call these $\textup{Hom}$\textit{-functors.}

Let $\mathcal{M}$ be the class of admissible monomorphisms in an exact category ($\catsa, \Ex$) and
$\mathcal{P}$ be the class of admissible epimorphisms in this category.
We shall write $\Ex=(\mathcal{M},\mathcal{P})$ for brevity. Then an object $I\in \mathcal{A}$ is
$\mathcal{M}$-injective precisely when $\Hom{-}{I}$ is exact as a functor to $(\Abcat, \Exmax)$ and
$P\in\mathcal{A}$ is $\mathcal{P}$-projective precisely when $\Hom{P}{-}$ is exact.

In an abelian category with enough injectives there is a standard notion of injective dimension of an object. This is defined using injective resolutions. In order for these methods to transfer to the setting of exact categories we need to take extra care in our set up.
To this end, we work with resolutions built from particular types of morphisms.
\begin{definition}\label{Definition: admissible morphism}
Let $(\mathcal{A},\Ex)$ be an exact category and $\Ex=(\mathcal{M},\mathcal{P})$. A morphism $f\in\Hom{E}{F}$ is called
\textit{admissible} if there exist an object $G\in\mathcal{A}$ and morphisms $\pi_f\in\classofmorphisms{P}{E}{G}$ and
$\mu_f\in\classMofmorphisms{G}{E}$ such that $f=\mu_f\pi_f$. (Such decomposition is unique up to isomorphism.)

A sequence of admissible morphisms
\[
\begin{tikzpicture}[auto] \matrix(m)[matrix of math nodes, column sep=2em, row sep=0.5em]{
{\cdots} &	{E_1} &\phantom{G_1}  &  {E_2} & \phantom{G_2} & {E_3} & \phantom{G_3} & {\cdots} \\
 {} &	\phantom{E_1} & {G_1}  &  \phantom{E_2} & {G_2} & \phantom{E_3} & {G_3} & \phantom{\cdots}   \\};
\draw[->] (m-1-1) to node {} (m-1-2); \draw[->] (m-1-2) to node {\footnotesize{$f_1$}} (m-1-4); \draw[->] (m-1-4) to node {\footnotesize{$f_2$}} (m-1-6); \draw[->] (m-1-6) to node {\footnotesize{$f_3$}} (m-1-8);
\draw[->>] (m-1-2) to node[swap] {\footnotesize{$\pi_{1}$}} (m-2-3); \draw[->>] (m-1-4)  to node[swap]  {\footnotesize{$\pi_{2}$}} (m-2-5); \draw[->>] (m-1-6) to node[swap]  {\footnotesize{$\pi_{2}$}} (m-2-7);
\draw[>->] (m-2-3) to node[swap] {\footnotesize{$\iota_{1}$}} (m-1-4); \draw[>->] (m-2-5)  to node[swap]  {\footnotesize{$\iota_{2}$}} (m-1-6); \draw[>->] (m-2-7) to node[swap]  {\footnotesize{$\iota_{3}$}} (m-1-8);
\end{tikzpicture}
\]
is called \textit{exact} if, for each $n\in\mathbb{Z},$ $(\iota_{n},\pi_{n+1})$ is in~$\Ex$.

An exact sequence of admissible morphisms of the form
\[
\begin{tikzpicture}[auto] \matrix(m)[matrix of math nodes, column sep=3em]{
	{E}\vphantom{I^0} & {I^0} &	{I^1} &   {I^2} & {\cdots}\vphantom{I^0} \\};
\draw[>->] (m-1-1) to node {\footnotesize{$\iota$}} (m-1-2); \draw[->] (m-1-2) to node {\footnotesize{$d^0$}} (m-1-3); \draw[->] (m-1-3) to node {\footnotesize{$d^1$}} (m-1-4); \draw[->] (m-1-4) to node {\footnotesize{$d^2$}} (m-1-5);
\end{tikzpicture}
\]
where each $I^n$ is $\mathcal{M}$-injective is called an \textit{$\mathcal{M}$-injective resolution of~$E$.}
\end{definition}

It is easy to see that, if $\mathcal{M}$ is the class of admissible monomorphisms in an exact category $(\mathcal{A},\Ex)$
and $\mathcal{A}$ has enough $\mathcal{M}$-injectives, then an injective resolution exists for every object in~$\mathcal{A}$.
In the same situation, a little more work gives the following analogue of the injective dimension theorem for modules over rings.

\begin{theorem}\label{Theorem: Injective dimension theorem}
Let $\mathcal{M}$ be the class of admissible monomorphisms in an exact category $(\mathcal{A},\Ex)$.
Suppose $\mathcal{A}$ has enough $\mathcal{M}$-injectives. The following are equivalent for $n\geq{1}$ and every $E\in\mathcal{A}$.
	\begin{enumerate}[label=\upshape(\roman*)]
	\item If there is an exact sequence of admissible morphisms
 \begin{equation}\label{Diagram: injective dimension theorem diagram 1}\begin{tikzpicture}[auto, baseline=(current  bounding  box.center)]
 		\matrix(m)[matrix of math nodes, column sep=3em]{
		{E}\vphantom{I^0} & {I^0} &	{\cdots}\vphantom{I^0} &   {I^{n-1}} & {F}\vphantom{I^0} \\};
	\draw[>->] (m-1-1) to node { } (m-1-2); \draw[->] (m-1-2) to node { } (m-1-3); \draw[->] (m-1-3) to node { } (m-1-4); \draw[->>] (m-1-4) to node { } (m-1-5);
	\end{tikzpicture}
	\end{equation}
	 with each $I^{m}$, $0\leq m\leq n-1$ injective, then $F$ must be injective;
	\item There is an exact sequence of admissible morphisms	
	\begin{equation}\label{Diagram: injective dimension theorem diagram 2}\begin{tikzpicture}[auto, baseline=(current  bounding  box.center)]
			\matrix(m)[matrix of math nodes, column sep=3em]{
		{E}\vphantom{I^0} & {I^0} &	{\cdots}\vphantom{I^0} &   {I^{n-1}} & {I^n}\vphantom{I^0} \\};
	\draw[>->] (m-1-1) to node { } (m-1-2); \draw[->] (m-1-2) to node { } (m-1-3); \draw[->] (m-1-3) to node { } (m-1-4); \draw[->>] (m-1-4) to node { } (m-1-5);
	\end{tikzpicture}
\end{equation}
with each $I^{m}$, $0\leq m\leq n$ injective.
\end{enumerate}	
\end{theorem}
\begin{remark}\label{Remark: About Injective dimension theorem}
That (i) implies (ii) in Theorem~\ref{Theorem: Injective dimension theorem} follows easily from the fact that every object in such a category $\mathcal{A}$ has an $\mathcal{M}$-injective resolution and by the definition of an exact sequence of admissible morphisms. To get (ii) implies (i), one can follow the method of the same result for general abelian categories, or indeed module categories, making slight adjustments when necessary. That is, one goes via the route of Ext groups. These are defined using derived functors arising from Hom-functors. A nice explanation of derived functors from exact categories can be found in \cite[Section~12]{B08}. The module category versions of the injective dimension theorem can be found in various places, see, for example, \cite[Lemma~4.1.8]{W94} or \cite[Proposition~4.8]{O2000}. An explicit proof of Theorem~\ref{Theorem: Injective dimension theorem} is given
in~\cite{Rosb2020}.
\end{remark}
We now introduce the homological dimension we are after.
\begin{definition}\label{Definition: Injective Dimension}
Let $E\in\mathcal{A}$ and $\mathcal{M}$ be the class of admissible monomorphisms in an exact category $(\mathcal{A},\Ex)$.
We say $E$ has \textit{finite $\mathcal{M}$-injective dimension} if there exists an exact sequence of admissible morphisms as in
Diagram~\eqref{Diagram: injective dimension theorem diagram 2} with all $I^m$ $\mathcal{M}$-injective. If $E$ is of finite
$\mathcal{M}$-injective dimension we write $\Midim{E}=0$ if $E$ is $\mathcal{M}$-injective and $\Midim{E}=n$ if $E$ is
not $\mathcal{M}$-injective and $n$ is the smallest natural number such that there exists an exact sequence of admissible morphisms as in Diagram~\eqref{Diagram: injective dimension theorem diagram 2} with all $I^m$ $\mathcal{M}$-injective.
If $E$ is not of finite $\mathcal{M}$-injective dimension, we write $\Midim{E}=\infty$.
 	
The \textit{global dimension} of the exact category $(\mathcal{A},\Ex)$ is
\[
\sup\setst{\Midim{E}}{E\in\mathcal{A}}\in\NN_0\cup\{\infty\}.
\]
\end{definition}
\begin{remark}
It is easy to observe, using Theorem~\ref{Theorem: Injective dimension theorem}, that the  $\mathcal{M}$-injective dimension of an
object $E$ in an exact category $(\mathcal{A},\Ex)$ is independent of any choice of injective resolutions.
That is, if there exists an exact sequence of admissible morphisms as in Diagram~\eqref{Diagram: injective dimension theorem diagram 2}
with all $I^m$ injective and also the following injective resolution of $E$
\begin{equation}\label{Diagram: injective resolution}
	\begin{tikzpicture}[auto, baseline=(current  bounding  box.center)]
		\matrix(m)[matrix of math nodes, column sep=3em]{
			{E}\vphantom{J^0} & {J^0} &	{\cdots}\vphantom{J^0} &   {J^{n-1}} & {J^n}\vphantom{J^0}  &{\cdots}\vphantom{J^0} 	\\};
		\draw[>->] (m-1-1) to node { } (m-1-2); \draw[->] (m-1-2) to node {  \footnotesize{$d^{0}$}} (m-1-3); \draw[->] (m-1-3) to node { } (m-1-4); \draw[->] (m-1-4) to node { \footnotesize{$d^{n-1}$}} (m-1-5); \draw[->] (m-1-5) to node { \footnotesize{$d^{n}$}} (m-1-6);
	\end{tikzpicture}
\end{equation}
then, by Theorem~\ref{Theorem: Injective dimension theorem}, $d^{n-1}$ must also be an admissible epimorphism.
\end{remark}

If $\mathcal{P}$ is the class of admissible epimorphisms then the dual procedure to the above yields the notion of projective resolutions and projective dimension for an object. There is a dual result to Theorem~\ref{Theorem: Injective dimension theorem} which leads to the notion of the projective dimension of an object. An object will have projective dimension $0$ if and only if it is $\mathcal{P}$-projective. We are able to characterise global dimension $0$ using projectives by Corollary~\ref{cor:all-injs-are-projs} below.

\begin{remark}\label{Remark: Global dimension using projective dimensions}
Let $\mathcal{M}$ be the class of admissible monomorphisms and $\mathcal{P}$ be the class of admissible epimorphisms in an exact category. If there are enough $\mathcal{P}$-projectives and enough $\mathcal{M}$-injectives, the supremum of all the $\mathcal{P}$-projective dimensions will coincide with the supremum of all the $\mathcal{M}$-injective dimensions. So we are always able to define the global dimension using only $\mathcal{P}$-projectives
provided there are enough of them.
\end{remark}

The following proposition is a classical result known as the `Splitting Lemma'.

\begin{prop}\label{Lemma: splitting lemma for additive categories}\textup{\textbf{(Splitting Lemma)}}
Let $\catsa$ be an additive category with kernels and cokernels and suppose
\[
\begin{tikzpicture}[auto]
		\matrix[ column sep={4em}]{	
			\node (E) {$E$}; & \node (F) {$F$};& \node (G) {$G$}; \\	};
		\draw[->] (E) to node {\footnotesize{$\mu$}} (F); 	\draw[->] (F) to node {\footnotesize{$\pi$}} (G);
	\end{tikzpicture}
\]
is a kernel-cokernel pair in~$\catsa$. Then the following are equivalent: 	
\begin{enumerate}[label=\upshape(\roman*)]
	\item The kernel-cokernel pair is split;
	 \item there exists  $\pi_E\in\Homunder{F}{E}{\mathcal{A}}$ such that $\pi_E\mu=\idntyof{E}$;
	 \item there exists  $\mu_F\in\Homunder{G}{F}{\mathcal{A}}$  such that $\pi\mu_F=\idntyof{G}$.
\end{enumerate}
\end{prop}
\begin{cor}\label{cor:all-injs-are-projs}
Let $(\mathcal{A},\Ex)$ be an exact category and $\Ex=(\mathcal{M},\mathcal{P})$.
Every object in $\mathcal{A}$ is $\mathcal{M}$-injective if and only if every object is $\mathcal{P}$-projective.
\end{cor}
\begin{proof}
This follows from the Splitting Lemma and Remark~\ref{Remark: Retractions of injectives and projectives}.
\end{proof}

\section{Operator Algebras with Global Dimension Zero}\label{Section: Relative Homological Algebra for Operator Modules over Operator Algebras}

\noindent
When $A$ is a unital operator algebra and $\mathcal{M}$ is the class of admissible monomorphisms in the exact category $(\OMODAinf,\Exmax)$,
it is unclear whether there are enough $\mathcal{M}$-injectives. The canonical object $\CB AI$, where $I$ is injective in $\Opinf$, which
is the analogue of the canonical injective object in algebraic module categories, lies in the larger category of matrix normed
modules \cite{Crann2017}, \cite{Rosb_arch2020}.
We will discuss this issue at the end of the present section.
However, there is a different exact structure we can place on $\OMODAinf$ that has enough projectives. In this section we discuss this
structure and give a characterisation of when the global dimension associated to this exact category is zero.

Let $\Exrel$ be the class of kernel-cokernel pairs in $\OMODAinf$ that split in~$\Opinf.$

\begin{prop}\label{Proposition: Erel forms an exact structure on OMODAinf}
Let $A$ be an operator algebra. Then $(\OMODAinf,\Exrel)$ is an exact category.
\end{prop}
\begin{proof}
By Theorem~\ref{Theorem: Emax in OMODAinf forms an exact structure}, $(\OMODAinf,\Exmax)$ is an exact category.
In particular, this holds for $\OMODinf{\mathbb{C}}=\Opinf$. The forgetful functor
$(\OMODAinf,\Exmax)\longrightarrow{(\Opinf,\Exmax)}$ is exact and as in Example~\ref{Example:Min class of kernel-cokernel pairs},
we know that $\Exmin$ forms an exact structure on~$\Opinf$.
The result follows from Proposition~\ref{Proposition: Using exact functors to get exact structures}.
\end{proof}
We next show that $(\OMODAinf,\Exrel)$ has enough projectives.
\begin{prop}\label{OMODAinf has enought P-projectives}
Let $A$ be a unital operator algebra. For every $E\in\Opinf\!$, the Haagerup tensor product $\hten EA$ is $\Prel$-projective.
In particular, $(\OMODAinf,\Exrel)$ has enough $\Prel$-projectives.
\end{prop}
\begin{proof}
Suppose we have the following diagram of morphisms in $\OMODAinf$,
\[
\begin{tikzpicture}[auto]
	\matrix(m)[row sep={1cm}, column sep=1cm, matrix of math nodes]
	 { { } & {\hten{E}{A}} \\ {F} & G\\ }; \draw[->] (m-1-2) to node {{\footnotesize{$f$}}} (m-2-2); \draw[->] (m-2-1) to node [swap]{{\footnotesize{$\pi$}}} (m-2-2);
\end{tikzpicture}
\]
with $\pi\in\Prel$ and that $\widetilde\pi\in\CB GF$ satisfies $\pi\widetilde\pi=\idntyof{G}$.
Let $\mapfromto{g'}{\hten EA}{F}$ be the map defined on elementary tensors by
$g'(x\otimes{a})=(\widetilde{\pi}f(x\otimes{1}))\cdot{a}$ for $x\in E, a\in A$.
Since $g'$ is the composition of
\[
\begin{tikzpicture}[auto]
\matrix(m)[row sep={1cm}, column sep=3.5em, matrix of math nodes]
{  {\hten{E}{A}} & \hten{(\hten{E}{A})}{A} & \hten{F}{A} & {F}\vphantom{_h}  \\ };
\draw[->] (m-1-1)to node {\footnotesize{$\iota\otimes{\idntyof{A}}$}} (m-1-2);
\draw[->] (m-1-2)to node {\footnotesize{$\widetilde{\pi}{f}\otimes{\idntyof{A}}$}} (m-1-3);
\draw[->] (m-1-3)to node {\footnotesize{$m$}} (m-1-4);
\end{tikzpicture}
\]
where $\iota\colon E\to\hten EA$ is $x\mapsto x\otimes1$ and $m$ is the completely contractive
linearisation of the module action of $A$ on~$F$, it is a well-defined completely bounded $A$-module map.
As
\[
\pi g'(x\otimes{a})=\pi(\widetilde{\pi}f(x\otimes{1})\cdot{a})=(\pi\widetilde{\pi}f(x\otimes{1}))\cdot{a}=f(x\otimes{1})\cdot{a}=f(x\otimes{a}).
\]
we have $\pi g'=f$.

Now let $E\in\OMODAinf$. Let $\mapfromto{P}{\hten{E}{A}}{E}$
be the canonical complete contraction given by $P(x\otimes{a})=x\cdot{a}$.
As each $x=[x_{i j}]\in M_n(E)$ is of the form $P_{n}([x_{i j}\otimes{1}])$ and $\norm{[x_{i j}\otimes{1}]}_{n}\leq\norm{[x_{ij}]}_n$,
$P$ is completely open.
For each $x\otimes{a}$ in $\hten{E}{A}$ and $a_0\in{A}$, we have
\[
P(x\otimes{a})\cdot{a_0}=(x\cdot{a})\cdot{a_0}=x\cdot{aa_0}=P(x\otimes{aa_0})=P((x\otimes{a})\cdot{a_0}),
\]
that is, $P$ is an $A$-module map. Define $\widetilde{P}\in\CB{E}{\hten{E}{A}}$ by ${\widetilde{P}(x)}=x\otimes{1}$.
Then $P\widetilde{P}=\idntyof{E}$, hence $P\in\Prel$ as required.
\end{proof}
Here comes the main result of the paper.
\begin{theorem}\label{Theorem: Relative Global Dimension 0 iff classically semisimple}
Let $A$ be a unital operator algebra. The following are equivalent:
\begin{enumerate}[label=\upshape(\roman*)]
\itemsep2pt
\item The global dimension of\/ $(\OMODAinf, \Exrel)$ is zero.
\item Every object in\/ $\OMODAinf$ is $\mathcal{M}_{rel}$-injective,
    where $\mathcal{M}_{rel}$ is the class of admissible monomorphisms in~$\Exrel$.
\item Every object in\/ $\OMODAinf$ is $\mathcal{P}_{rel}$-projective,
    where $\mathcal{P}_{rel}$ is the class of admissible epimorphisms in $\Exrel$.
\item $A$ is classically semisimple.
\end{enumerate}
\end{theorem}
In keeping with traditional notation we shall write $\dgcb A$ for the global dimension of
$(\OMODAinf, \Exrel)$ and call it \textit{the completely bounded global dimension\/} of the operator algebra~$A$.
\begin{cor}
For every unital operator algebra~$A$, its completely bounded global dimension $\dgcb A$ is zero if and only if $A$
is a finite direct sum of full matrix algebras.
\end{cor}
The equivalence of (i) and (ii) in Theorem~\ref{Theorem: Relative Global Dimension 0 iff classically semisimple}
follows from Definition~\ref{Definition: Injective Dimension} and (ii)${}\iff{}$(iii) is a
consequence of the splitting lemma (Corollary~\ref{cor:all-injs-are-projs}).

The argument for the implication (iv) ${}\Rightarrow{}$ (iii) follows closely the proof of
\cite[Theorem 3.5.17]{HELhomBook}; for the details see \cite{Rosb2020}, Section~5.3.
A direct argument of (iv) ${}\Rightarrow{}$ (ii) is provided by Proposition~\ref{Prop: Semisimple implies all rel injective} below.
The final implication (iii) ${}\Rightarrow{}$ (iv) needs some preparation.

Our approach is based on the following well-known characterisation; see, e.g., \cite[Theorem~4.40]{O2000}.
\begin{prop}\label{Useful characterisation of semisimplicity of a ring}
A unital algebra is classically semisimple if and only if each of its maximal submodules is a direct summand.
\end{prop}
In other words, to show that a unital algebra $A$ is classically semisimple, we have to show that,
for every maximal right ideal~$I$ of~$A$, the exact sequence of right $A$-modules
\begin{equation}\label{eq:ss-max-right-ideals}
\begin{tikzpicture}[auto]\matrix(m)[matrix of math nodes, column sep=1.2cm]
{0\vphantom{A/I} & I\vphantom{A/I} & A\vphantom{A/I} & A/I & 0\vphantom{A/I} \\};
\draw[->] (m-1-1) to node {} (m-1-2); \draw[->](m-1-2) to node {{\footnotesize{$\iota$}}} (m-1-3); \draw[->](m-1-3) to node {{\footnotesize{$\pi$}}} (m-1-4); \draw[->](m-1-4) to node {} (m-1-5);
\end{tikzpicture}
\end{equation}
splits. As $A/I$ is a simple right $A$-module, this is evidently equivalent to the existence of an $A$-module map
$f\colon A/I\to A$ such that $f(A/I)\nsubseteq I$ (cf.\ Lemma~\ref{Lemma: splitting lemma for additive categories}).

Let $A$ be a unital operator algebra; as every maximal right ideal $I$ of $A$ is closed, we can consider
the sequence~\eqref{eq:ss-max-right-ideals} above within $\OMODAinf$. The next auxiliary result will enable us to complete
the final step in the proof in the proof of Theorem~\ref{Theorem: Relative Global Dimension 0 iff classically semisimple}.
\begin{lem}\label{lem:form of r_I and fact its twosided ideal}
Let $I$ be a closed right ideal in the unital operator algebra~$A$.
Set
\begin{equation*}
\begin{split}
S_1  &=\setst{f(x)}{x\in A/I, \mapfromto{f}{A/I}{A} \text{ is an }A\text{-module map} },\\
S_2  &=\setst{ab}{a,b\in A,\ ay=0 \text{ for all } y\in I}.
\end{split}
\end{equation*}
Then  $S_1=S_2$ and its linear span $S$ is a two-sided ideal.
\end{lem}
\begin{proof}
That the linear span of $S_2$ is a two-sided ideal is obvious.
Let $\mapfromto{\pi}{A}{A/I}$ be the canonical quotient map.
Let $a\in A$. We claim there exists an $A$-module map $\mapfromto{f}{A/I}{A}$ such that
$f(\pi(1_A))=a$ if and only if $ay=0$ for all $y\in I$.
		
Suppose there exists such an $f$. Take $y\in I$.
Then $ay=f(\pi(1_A))y=f(\pi(y))=f(0)=0$. Now suppose $ay=0$ for all $y\in I$.
Define $\mapfromto{f}{A/I}{A}$ by $f(x)=ab$, where $\pi(b)=x$.
Clearly this is well defined and $f(\pi(1_A))=a$. So the claim is true.
The result then follows as an element is of the form $f(x)$ if and only if it is of the
form $f(\pi(1_A)){b}$ for $b\in A$ such that $\pi(b)=x$.
\end{proof}
The proof of the next result mirrors \cite[Proposition IV.4.4.]{HELhomBook}.
\begin{prop}\label{using hten and rel projectivity}
Let $\mapfromto{\phi}{A}{B}$ be a unital completely contractive homomorphism between unital operator algebras
and let $E\in\OMODAinf$ be $\mathcal{P}_{rel}$-projective. Suppose $F\in\OMODinf{B}$ and
$\mapfromto{f}{E}{F}$ is a completely bounded $A$-module map, where we consider $F$ as an $A$-module by the restriction of scalars via~$\phi$.
For any $x\in E$ with $f(x)\neq0$, there exists a completely bounded $A$-module map
$\mapfromto{g}{E}{A}$ such that $\phi{g}(x)\neq0$.
\end{prop}
\begin{proof}
In the following proof, for any $G\in\OMODinf{C}$, where $C$ is an operator algebra,
we will denote by $P_{G,C}$ the completely contractive $C$-module map $\hten{G}{C}\rightarrow{G}$ defined on elementary tensors by $P_{G,C}(z\otimes{c})=z\cdot{c}$ (where $z\in G, c\in C$). If $G$ is a left operator $C$-module, we denote the similarly defined completely contractive $C$-module map by $\mapfromto{{}_{C,G}P}{\hten{C}{G}}{G}$.
		
By our assumptions, we have the following commutative diagram of completely bounded linear maps.
		\[\begin{tikzpicture}[auto]\matrix(m)[matrix of math nodes, column sep=2cm, row sep=1cm] {E{_{\vphantom{h}}} & \hten{E}{A} \\  F{_{\vphantom{h}}} & \hten{F}{B} \\ };
		\draw[->] (m-1-2) to node {{\footnotesize{$P_{E,A}$}}} (m-1-1);  	\draw[->] (m-1-1) to node [swap]{{\footnotesize{$f$}}} (m-2-1);
		\draw[->] (m-2-2) to node {{\footnotesize{$P_{F,B}$}}} (m-2-1);    \draw[->] (m-1-2) to node {{\footnotesize{$f\otimes{\phi}$}}} (m-2-2);
		\end{tikzpicture}\]
As in the proof of Proposition~\ref{OMODAinf has enought P-projectives}, we have $P_{E,A}\in\mathcal{P}_{rel}$.
Hence, there exists a completely bounded $A$-module map $\mapfromto{\iota}{E}{\hten{E}{A}}$ such that $P_{E,A}\,\iota=\idntyof{E}$.
Then it is clear that the above diagram stays commutative if we replace $P_{E,A}$ with~$\iota$.
Moreover, for any (completely) bounded linear functional $\mapfromto{\alpha}{F}{\mathbb{C}}$ it is easy to see the following
diagram of completely bounded linear maps is commutative:
		\[\begin{tikzpicture}[auto]
		\matrix(m)[matrix of math nodes, column sep=2cm, row sep=1cm]
		{E{_{\!\!\phantom{h}}} & \hten{E}{A} & \hten{\mathbb{C}}{A} & {A{_{\!\!\phantom{h}}}} \\
			F{_{\!\!\phantom{h}}} & \hten{F}{B} & \hten{\mathbb{C}}{B} & {B{_{\!\!\phantom{h}}}} \\ };
		
		\draw[->] (m-1-1) to node {{\footnotesize{$\iota$}}} (m-1-2);
		\draw[->] (m-1-2) to node {{\footnotesize{$\alpha f\otimes\idntyof{A}$}}} (m-1-3);
		\draw[->] (m-1-3) to node {{\footnotesize{${}_{\mathbb{C}, A}P$}}} (m-1-4);
		
		\draw[->] (m-1-1) to node [swap]{{\footnotesize{$f$}}} (m-2-1);
		\draw[->] (m-1-2) to node {{\footnotesize{$f\otimes{\phi}$}}} (m-2-2);
		\draw[->] (m-1-3) to node [swap]{{\footnotesize{$\idntyof{\mathbb{C}}\otimes{\phi}$}}} (m-2-3);
		\draw[->] (m-1-4) to node [swap]{{\footnotesize{$\phi$}}} (m-2-4);

		\draw[->] (m-2-2) to node {{\footnotesize{$P_{F,B}$}}} (m-2-1);
		\draw[->] (m-2-2) to node [swap]{{\footnotesize{$\alpha\otimes\idntyof{B}$}}} (m-2-3);
		\draw[->] (m-2-3) to node [swap]{{\footnotesize{${}_{\mathbb{C}, B}P$}}} (m-2-4);
		\end{tikzpicture}\]
Let $a\in A$ and $x\otimes{a'}$ be an elementary tensor in $\hten{E}{A}$ and $\lambda\otimes{a''}$ be an elementary tensor in $\hten{\mathbb{C}}{A}$.
Then
\[
((\alpha f\otimes\idntyof{A})(x\otimes{a'}))\cdot{a}=\alpha f(x)\otimes{a'a}=(\alpha f\otimes\idntyof{A})((x\otimes{a'})\cdot{a})
\]
and $_{\mathbb{C}, A}P(\lambda\otimes{a''})\cdot{a}=\lambda{a''a}= {}_{\mathbb{C}, A}P((\lambda\otimes{a''})\cdot{a})$.
By continuity and linearity, $\alpha f\otimes\idntyof{A}$ and ${}_{\mathbb{C}, A}P$ are $A$-module maps. So, for any linear functional $\alpha$,
we have that $\mapfromto{{}_{\mathbb{C}, A}P(\alpha f\otimes\idntyof{A})\iota}{E}{A}$ is a completely bounded $A$-module map.
We now show that there exists an $\alpha$ that makes this the desired $A$-module map.
		
Let $x\in E$ such that $f(x)\neq0$.
Then $P_{F,B}(f\otimes\phi)\iota(x)\neq0$; in particular $u=(f\otimes\phi)\iota(x)\in\hten{F}{B}$ is non-zero.
By Lemma~\ref{lem:isom-embeddings}, there exist $\alpha\in F^*$, $\beta\in B^*$ such that $(\alpha\otimes\beta)(u)\neq0$.
As $\alpha\otimes\beta$ is the composition $(\idntyof{E}\otimes\beta)(\alpha\otimes\idntyof{B})$  we obtain that
$(\alpha\otimes\idntyof{B})u\neq0$ and therefore ${}_{\mathbb{C}, B}P(\alpha\otimes\idntyof{B})u\neq0$.
		
Let $g={}_{\mathbb{C}, A}P(\alpha(f)\otimes\idntyof{A})\iota$. Commutativity of the diagram gives then that $\phi{g}(x)\neq0$ as required.	
\end{proof}
We can now complete the \textit{Proof\/} {of Theorem~\ref{Theorem: Relative Global Dimension 0 iff classically semisimple}
{\rm(iii)${}\Rightarrow{}$(iv)}:}

\smallskip\noindent
Let $I$ be a maximal right ideal of~$A$.
Suppose that the image of every $A$-module map $A/I\rightarrow{A}$ is a subset of $I$; then, with notation
as in Lemma~\ref{lem:form of r_I and fact its twosided ideal}, $S_1\subseteq{I}$, and hence $S\subseteq{I}$.
Put $E=A/I$ and $P=\setst{a\in A}{x\cdot{a}=0, \forall x\in E},$ the right annihilator of~$E$.
Then $P$ is a closed two-sided ideal in~$A$. Moreover, $S\subseteq P$. Indeed, suppose $a\in S$ and $x\in E$ is of the form $\pi(b)$ for some $b\in A$, where $\mapfromto{\pi}{A}{E}$ is the canonical epimorphism.
Then $x\cdot{a}=\pi(b)\cdot{a}=\pi(ba)=0$ because $S$ is a left ideal contained in~$I$.

Put $B=A/P$ and let $\mapfromto{\phi}{A}{B}$ be the canonical epimorphism of unital operator algebras.
Then $E$ is a right $B$-module with action defined in the following way:
for $x\in E$ and $b\in B$, $x\cdot{b}=\pi(aa')$ where $x=\pi(a)$, $b=\phi(a')$.
This is well defined by the definition of~$P$.
Then $x\cdot{b}=x\cdot{a'}$.
In fact, $E\in\OMODinf{B}$; suppose $x\in M_n(E)$, $b\in M_n(B)$. For every $a'\in M_n(A)$ with $b=\phi_n(a')$, $\norm{x\cdot{b}}_n=\norm{x\cdot{a'}}_n\leq\norm{x}_n\norm{a'}_n,$ so $\norm{x\cdot{b}}_n\leq\norm{x}_n\norm{b}_n$.
Let $F=E$ and $f=\idntyof{E}$; then obviously we have that $f(\pi(1_A))\neq0$.
If $E$ is $\mathcal{P}_{rel}$-projective, by Proposition~\ref{using hten and rel projectivity}, for each $x\notin I$,
there exists an $A$-module map $\mapfromto{g}{E}{A}$ such that $\phi(g(\pi(x)))\neq0$, so $g(\pi(x))\notin{P}$.

Since, by (the proof of) Lemma~\ref{lem:form of r_I and fact its twosided ideal},
$g(\pi(x))\in S\subseteq P$ it follows that $E$ cannot be $\mathcal{P}_{rel}$-projective.
Consequently, if $A$ is not classically semisimple, using Proposition~\ref{Useful characterisation of semisimplicity of a ring},
we conclude that not all modules in $\OMODAinf$ can be $\mathcal{P}_{rel}$-projective.
That is, if (iv) does not hold, then (iii) cannot hold either.\qed

\begin{remark}\label{Remark:conclusions}
Implication (iii)${}\Rightarrow{}$(iv) in Theorem~\ref{Theorem: Relative Global Dimension 0 iff classically semisimple}
answers Helemskii's question, for operator algebras, in \cite[Section~7]{Helem2009} in the positive since it is easy to see that
his relative structure is equivalent to ours.
Relative homological algebra is common in the ring theory setting, cf., e.g., \cite[Chapter~IX]{HiltStamm}
or \cite[Section V.7]{Mitch}.
Paulsen undertook a systematic study in the setting of operator modules in~\cite{Paulsen1998}, see also~\cite{FrankPaulsen2003}.
He discovered an intimate interrelation between the cohomology groups that arise in this relative theory (which is equivalently described
here by our relative exact structure) and the completely bounded Hochschild cohomology groups.
See in particular Propositions~5.5 and~6.4 in~\cite{Paulsen1998}.
Specialising to bimodules it then follows that all operator $A$-bimodules over a unital \C* $A$
are relatively projective if and only if $A$ is finite dimensional which in turn is equivalent to $A$ possessing a diagonal
\cite[Theorem~6.13 and Corollary~6.14]{Paulsen1998}. For an extension of the latter to unital operator algebras, see~\cite{PauSm2002}.
\end{remark}
In the remainder of this section we shall discuss the interrelations between the various types of injectivity within
$\OMODAinf$ but also in comparison to the larger category $\mnMODAinf$. It will become clear that injectivity is not
determined by the category but rather the exact structure which one puts on the category.
For a discussion of injectivity in general terms, we refer to~\cite{MR20}.

Let $\XMODinfA$ be any additive category whose objects are operator spaces which are right $A$-modules and whose morphisms are the completely bounded $A$-module maps between these objects, where $A$ is a unital operator algebra (in the sense of \cite[Section~2.1]{BLM04}).
Suppose $\XMODinfA$ is closed under direct sums, closed submodules and quotients.
Then the kernels (respectively, cokernels) in $\XMODinfA$ are as described in Proposition~\ref{Proposition: Description of kernels in OMODAinf} (resp., Proposition~\ref{Proposition: Description of cokernels in OMODAinf}).
Moreover, the proofs of Propositions~\ref{Proposition: Pushouts and kernels in OMODAinf} and Proposition~\ref{Proposition: Pullbacks and cokernels in OMODAinf} still work for $\XMODinfA$
and therefore, $(\XMODinfA, \Exmax)$ and $(\XMODinfA, \Exrel)$ are exact categories,
where $\Exrel$ again denotes the kernel-cokernel pairs that split in~$\Opinf$ (Proposition~\ref{Proposition: Erel forms an exact structure on OMODAinf}).

In particular, this holds for the category $\mnMODAinf$, whose objects are the non-degenerate matrix normed $A$-modules.
An operator space $E$ that is also a right $A$-module is known as a \textit{matrix normed $A$-module} if
the $A$-module action induces a completely contractive linear mapping $\projten{E}{A}\rightarrow{E},$
where $\otimesfrown$ denotes the operator space projective tensor product.
Comparing $\otimesfrown$ with $\otimes_h$ immediately tells us that $\OMODAinf$ is a full subcategory of $\mnMODAinf.$
(see \cite[Example 3.1.5]{BLM04}). The category $\mnMODAinf$ is, e.g., used in \cite{Aristov2002}, \cite{Crann2017} and \cite{Rosb_arch2020}.

For each $E\in\Opinf$, on $\CB{A}{E}$  we define the right $A$-module action by
$(T\cdot{a})(b)=T(ab),$ for all $T\in\CB{A}{E}$ and  $a,b\in A$.
Then $\CB{A}{E}$ is an object in $\mnMODAinf$ which is just a specific case of \cite[3.5.2]{BLM04}.
\begin{prop}\label{Proposition: Typical relinjs in mnMODainf}
Let $A$ be a unital operator algebra and let $\Nrel$ be the class of admissible monomorphisms in $(\mnMODAinf,\Exrel).$
For any $G\in \Opinf$, the matrix normed module $\CB{A}{G}$ is $\Nrel$-injective in $\mnMODAinf$.
Moreover, $E\in\mnMODAinf$ is $\Nrel$-injective if and only if it is a retract of some $\CB{A}{G}$.
\end{prop}
\begin{proof}
Suppose $E,F\in\mnMODAinf, G\in\Opinf$. Let
$\mu\in\morphismsnomathcal{\Nrel}{E}{F}$ and $f\in\CBmod{A}{E}{\CB{A}{G}}$.
We show that $\CB{A}{G}$ is $\Nrel$-injective by finding a morphism $g\in\CBmod{A}{F}{\CB{A}{G}}$
such that $f=g\mu$.

As $\mu\in\Nrel$, there exists $\widetilde{\mu}\in\CB{F}{E}$ such that ${\widetilde{\mu}}\mu=\idntyof{E}$.
For $y\in F$ write $f(\widetilde{\mu}(y\cdot{a}))(1_A)$ as $g(y)(a)$, for each $a\in A$.
This defines a completely bounded linear map $g(y)\in\CB{A}{G}$.

It is routine to verify that, in fact, this yields a morphism
$g\in\CBmod{A}{F}{\CB{A}{G}}$ and moreover, for all  $x\in E, a\in A$
\[\begin{split}
g(\mu(x))(a)=f(\widetilde{\mu}(\mu(x)\cdot{a}))(1)=
f(\widetilde{\mu}(\mu(x\cdot{a}))(1) &=f(x\cdot{a})(1)
\\
&=(f(x)\cdot{a})(1)=f(x)(a).
\end{split}\]
So $g\mu=f$ as required.

Now suppose that $E\in\mnMODAinf$ and define $\mapfromto{\iota}{E}{\CB{A}{E}}$ by
$\iota(x)(a)=x\cdot{a}$ for each $x\in E, a\in A$.
Clearly, $\iota$ is a completely isometric $A$-module map and thus a kernel map in $\mnMODAinf$.
Define $\mapfromto{\widetilde{\iota}}{\CB{A}{E}}{E}$ by
$\widetilde{\iota}(T)= T(1_A)$ for all $T\in\CB{A}{E}$.
Then $\widetilde{\iota}$ is a completely bounded linear map such that
${\widetilde{\iota}}\iota=\idntyof{E}$.
That is, $\iota\in\Nrel$.
The result follows by Remark~\ref{Remark: Retractions of injectives and projectives}.
\end{proof}

The above proposition is obtained in \cite[Section 2]{Crann2017} in a similar way but without making the categorical setting explicit.
Since we shall compare injectivity in different categories, we need to make sure our arguments fit the correct situation.

\begin{remark}\label{Remark: Rel injectivity in matrix normed case}
Let $\Mrel$ denote the class of admissible monomorphisms in $(\OMODAinf, \Exrel)$
and $\Nrel$ be the admissible monomorphisms in $(\mnMODAinf,\Exrel)$.
If $E,F,I\in\OMODAinf$, then $\morphismsnomathcal{\Mrel}{E}{F}=\morphismsnomathcal{\Nrel}{E}{F}$
and
\[
\Morunder{F}{I}{\OMODAinf}=\CBmod{A}{F}{I}=\Morunder{F}{I}{\mnMODAinf}.
\]
Hence, if $I$ is $\Nrel$-injective it must also be $\Mrel$-injective.
Therefore, if $E\in\OMODAinf$ and $E$ is a retract (in $\mnMODAinf$) of $\CB{A}{E}$, then $E$ will be $\Mrel$-injective.
\end{remark}
The next result shows that `completely bounded global dimension zero' does not depend on which of the two categories
one chooses.
\begin{prop}\label{Prop: Semisimple implies all rel injective}
Suppose $A$ is a classically semisimple unital operator algebra.
Then every object in $\OMODAinf$ is $\Mrel$-injective and every object in $\mnMODAinf$ is $\Nrel$-injective.
\end{prop}
\begin{proof}
By Proposition~\ref{Proposition: Typical relinjs in mnMODainf} and Remark~\ref{Remark: Rel injectivity in matrix normed case},
it suffices to show that for every $E\in\mnMODAinf,$
there exist $r\in\CBmod{A}{\CB{A}{E}}{E}$ and $s\in \CBmod{A}{E}{\CB{A}{E}}$
such that $rs=\idntyof{E}.$
First we will fix some notation:
		
There exists $n\in\mathbb{N}$ such that
$A=M_{m_1}(\mathbb{\mathbb{C}})\oplus M_{m_2}(\mathbb{\mathbb{C}})\oplus \cdots \oplus M_{m_n}(\mathbb{\mathbb{C}})$.
For each $k\in\{1,\ldots, n\}$ and $i,j\in\{1,\ldots, m_k\}$, let $e^k_{ij}$ denote the $n$-tuple in $A$ with all zero entries apart from the $k$-th entry which is a matrix in $M_{m_k}(\mathbb{C})$ with $1$ for the $ij$-th entry and $0$ everywhere else.
		
Note that $A$ is the linear span of the elements $e^k_{ij}$.
Moreover, $\norm{e^k_{ij}}=1$ and $e^k_{ij}e^{\ell}_{pq}=0$
unless $j=p,$ and $k={\ell}$ in which case $e^k_{ij}e^{\ell}_{pq}=e^k_{iq}$.
Then  $1_A=\sum_{k=1}^{n}\sum_{i=1}^{m_{k}}e^k_{ii}$ and for any $x\in E$ we have
$x=\sum_{k=1}^{n}\sum_{i=1}^{m_{k}}x\cdot{e^k_{ii}}$.
		
Let $\mapfromto{s}{E}{\CB{A}{E}}$ be the completely bounded $A$-module map defined by $s(x)(a)=x\cdot{a}$ for all $x\in E$, $a\in A$.
For each $T\in \CB{A}{E}$, let
\[
	r(T)= \sum_{k=1}^{n}\sum_{i=1}^{m_{k}}T(e^k_{i1})\cdot e^k_{1i}.
\]
It is clear that this defines a linear mapping $\mapfromto{r}{\CB{A}{E}}{E}$
and so will be an $A$-module map if for each $T\in \CB{A}{E}$
we have $r(T)\cdot{e^{\ell}_{pq}}=r(T\cdot{e^{\ell}_{pq}})$
for arbitrary $t\in\{1,\ldots, n \}$ and $p,q\in \{1,\ldots, m_k\}$.
		
We compute the two terms in question:
\begin{equation}\label{eqn: r(T.a)}
	r(T\cdot{e^{\ell}_{pq}})=\sum_{k=1}^{n}\sum_{i=1}^{m_{k}}T(e^{\ell}_{pq}e^k_{i1})\cdot e^k_{1i}=T(e^{\ell}_{pq}e^{\ell}_{q1})\cdot{e}^{\ell}_{1q}=T(e^{\ell}_{p1})\cdot {e}^{\ell}_{1q},
\end{equation}
as every other term is zero.
Similarly, as $e^k_{1i}e^{\ell}_{pq}=0$ unless $k={\ell}$ and $i=p,$ we have
\begin{equation}\label{eqn: r(T).a}
r(T)\cdot{e^{\ell}_{pq}}=\sum_{k=1}^{n}\sum_{i=1}^{m_{k}}T(e^k_{i1})\cdot e^k_{1i}e^{\ell}_{pq}= T(e^{\ell}_{p1}) \cdot (e^{\ell}_{1p}e^{\ell}_{pq})
	=T(e^{\ell}_{p1})\cdot {e}^{\ell}_{1q}.
\end{equation}
Comparing equations (\ref{eqn: r(T.a)}) and (\ref{eqn: r(T).a}) gives us that $r$ is an $A$-module map.
		
Let $x\in E$. Then
\[
 rs(x)= \sum_{k=1}^{n}\sum_{i=1}^{m_{k}}s(x)(e^k_{i1})\cdot e^k_{1i}=\sum_{k=1}^{n}\sum_{i=1}^{m_{k}}(x\cdot e^k_{i1})\cdot e^k_{1i} =\sum_{k=1}^{n}\sum_{i=1}^{m_{k}}x\cdot{e^k_{ii}} =x,
\]
so $rs=\idntyof{E}$ and all that remains is to show $r$ is completely bounded.
		
Note that for $T\in\CB{A}{E},$
\[
 r(T)=\sum_{k=1}^{n}\sum_{i=1}^{m_{k}} r^k_{i}(T),
\]
where $r^k_{i}(T)$ is defined to be $T(e^k_{i1})\cdot e^k_{1i}$.
Each $\mapfromto{r^k_{i}}{\CB{A}{E}}{E}$ is a linear map. Hence, it suffices that each $r^k_{i}$ is completely bounded.
		
Let $k\in\{1,\ldots, n \}$ and $i\in \{1,\ldots, m_k\}$.
For each $N\in\mathbb{N}$, let $e_N\in M_N(A)$ be the matrix with $e^k_{1i}$ as every entry in the leading diagonal
and zero everywhere else.
Note $\norm{e_N}_N=\norm{e^k_{1i}}=1$. For each $T=[T_{vw}]\in M_N(\CB{A}{E})$ we have
\[
 \bignorm{[(r^k_i)_N(T_{vw})]}_N=\bignorm{[T_{vw}(e^k_{i1})\cdot e^k_{1i}]}_N\leq\bignorm{[T_{vw}(e^k_{i1})]}_N\,\norm{e_N}_N
    \leq \norm{T}_N\leq\cbnorm{T}.
\]
So $r^k_i$ is completely contractive and $r$ is completely bounded as required.
\end{proof}
\begin{cor}\label{Corollary: Equivalence of dimension 0}
Let $A$ be a unital operator algebra.
The following are equivalent.
\begin{enumerate}[label=\upshape(\roman*)]
	\itemsep2pt
	\item $\dgcb A=0$.
	\item $A$ is classically semisimple.
	\item The global dimension of\/ $(\mnMODAinf, \Exrel)$ is zero.
\end{enumerate}
\begin{proof}
Statements (i) and (ii) are equivalent by Theorem~\ref{Theorem: Relative Global Dimension 0 iff classically semisimple}.
That (ii) implies (iii) is Proposition~\ref{Prop: Semisimple implies all rel injective}.
Finally, we have the implication (iii)${}\Rightarrow{}$(i) by Remark~\ref{Remark: Rel injectivity in matrix normed case}.
\end{proof}
\end{cor}
It is interesting to compare the relation between `global injectivity', that is, injectivity with respect to the maximal
exact structure and relative injectivity. To this end we record the following general result which is part of the
`injective version' of \cite[Proposition~11.3]{B08} where it is obtained for projective objects.
\begin{lem}\label{lem:inj-iff-abs-retract}
Let $\mathcal{(A,\Ex)}$ be an exact category and $\mathcal{M}$ be the class of admissible monomorphisms.
Then $I\in\mathcal{A}$ is $\mathcal{M}$-injective if and only if it is an absolute $\mathcal{M}$-retract,
that is, for every $\mu\in\classMofmorphisms{E}{F}$ with $F\in\mathcal{A}$ there exists ${\nu}\in\Mor{F}{E}$
such that $\nu\mu=\idntyof{E}$.
\end{lem}
\begin{proof}
That every $\mathcal{M}$-injective is an absolute retract is immediate from the definition.
Suppose $I\in\mathcal{A}$ is an absolute $\mathcal{M}$-retract and that $E,F\in\mathcal{A}$ with morphisms
$\mu\in\classMofmorphisms{E}{F}$ and $f\in\Mor{E}{I}$ are given.
By axiom [E2] of Definition~\ref{Definition: Exact Category}, there exists a (commutative) pushout square:
\[\begin{tikzpicture}[auto]
\matrix[matrix of math nodes, column sep=4em, row sep=3em](m)
{
	E & F \\
	I & C\\
};
\draw[>->] (m-1-1) to node {\footnotesize{$\mu$}} (m-1-2);
\draw[>->] (m-2-1) to node [swap]{\footnotesize{$\mu'$}} (m-2-2);
\draw[->] (m-1-1) to node [swap]{\footnotesize{$f$}} (m-2-1);
\draw[->] (m-1-2) to node {\footnotesize{$g$}} (m-2-2);
\draw[->, dashed] (m-1-2) to node {\footnotesize{$\nu g$}} (m-2-1);
\end{tikzpicture}
\]
such that $\mu'\in\classMofmorphisms{I}{C}.$
As $I$ is an absolute $\mathcal{M}$-retract, there exists $\nu\in\Mor{C}{I}$
such that $\nu\mu'=\idntyof{I}$.
Then $\nu{g}\in\Mor{F}{I}$ with $(\nu{g})\mu= \nu\mu'f=f$ as required.
\end{proof}
The next result can also be obtained in the general setting; however this would require to fix two categories
and four exact structures. In order to avoid such unnecessary generality, we restrict our attention to the two
categories of operator space modules, $\OMODAinf$ and $\mnMODAinf$.
We will suppress the formal forgetful functor from either of these categories to $\Opinf$ and simply consider
the modules as operator spaces when needed. But it is essential to note that the morphisms in $\Opinf$ are the completely bounded
linear mappings \textit{and not the complete contractions}. Therefore $\Opinf$-injectivity is \textit{not\/}
what one usually calls `injective operator space'.

Let $\catsa$ denote either $\OMODAinf$ or $\mnMODAinf$, for a unital operator algebra~$A$.
Since the morphisms, the admissible monomorphisms and the admissible monomorphisms in the relative structure are the same in
both categories, we can simply talk about `injective object' and `relatively injective object' in $\catsa$ below.
\begin{prop}\label{Prop: relative and N inj imply injective}
Let $E\in\mathcal{A}$. Suppose $E$ is relatively injective and $\Opinf$-injective. Then $E$ is injective.
\end{prop}
\begin{proof}
Take $\mu\in\classMofmorphisms{E}{F}$ for some $F\in\catsa$. By the lemma above,
we are done if there exists $\nu\in\CBmod AFE$ such that $\nu\mu=\idntyof{E}$.
As $\mu$ is an admissible monomorphism there exists a kernel-cokernel pair in $\Exmax$:
\begin{equation}\label{eq:ker-cok-global}
\begin{tikzpicture}[auto, baseline=(current  bounding  box.center)]
\matrix[matrix of math nodes, column sep=2.5em,](m)
{
	E & F & G\\
};
\draw[>->] (m-1-1) to node {\footnotesize{$\mu$}} (m-1-2);
\draw[->>] (m-1-2) to node {\footnotesize{$\pi$}} (m-1-3);
\end{tikzpicture}
\end{equation}
which gives a kernel-cokernel pair in $(\Opinf,\Exmin)$
\begin{equation}\label{eq:ker-cok-relative}
\begin{tikzpicture}[auto, baseline=(current  bounding  box.center)]
\matrix[matrix of math nodes, column sep=2.5em,](m)
{
	E & F & {G.}\\
};
\draw[>->] (m-1-1) to node {\footnotesize{$\mu$}} (m-1-2);
\draw[->>] (m-1-2) to node {\footnotesize{$\pi$}} (m-1-3);
\end{tikzpicture}
\end{equation}
Since $E$ is $\Opinf$-injective, $E$ is an absolute retract in $\Opinf$ (Lemma~\ref{lem:inj-iff-abs-retract}); hence there exists
$\theta\in\CB FE$ such that $\theta\mu=\idntyof{E}$.
By the Splitting Lemma (Lemma~\ref{Lemma: splitting lemma for additive categories}),
$(\mu,\pi)\in\Exmin$ in $\Opinf$, in \eqref{eq:ker-cok-relative}and therefore $(\mu,\pi)\in\Ex_{rel}$ in \eqref{eq:ker-cok-global}.
In particular $\mu\in\Mrel$. As $E$ is relatively injective,
by the other implication in Lemma~\ref{lem:inj-iff-abs-retract}, there exists $\nu\in\CBmod AFE$
such that $\nu\mu=\idntyof{E}$ as required so that $E$ is injective in~$\catsa$.
\end{proof}
With the same caveats as above we obtain a converse under an additional assumption.
\begin{prop}\label{Prop: injective implies others condition}
Suppose that, for every $E\in\catsa$, there exists $\mu\in\CBmod{A}{E}{J_E}$
for some $J_E\in\catsa$ such that $J_E$ is $\Opinf$-injective.
If $I\in\mathcal{A}$ is injective then $I$ is relatively injective as well as $\Opinf$-injective.
\end{prop}
\begin{proof}
Let $I\in\mathcal{A}$ be injective in~$\catsa$; then it is clearly relatively injective.
By assumption, there exist $J_I\in\mathcal{A}$ which is $\Opinf$-injective and
a kernel cokernel pair in $\Exmax$:
\[
	 \begin{tikzpicture}[auto]
	 \matrix[matrix of math nodes, column sep=2.5em,](m)
	 {
	 	{I}\vphantom{_I} & {J_I} & G\vphantom{_I}\\
	 };
	 \draw[>->] (m-1-1) to node {\footnotesize{$\mu$}} (m-1-2);
	 \draw[->>] (m-1-2) to node {\footnotesize{$\pi$}} (m-1-3);
	 \end{tikzpicture}
\]
By Lemma~\ref{lem:inj-iff-abs-retract}, there exists ${\nu}\in\CBmod{A}{J_I}{I}$ such that $\nu\mu=\idntyof{I}$.
This identity persists in $\Opinf$ so that $I$ is a retract of the $\Opinf$-injective operator space~$J_I$.
Hence $I$ is injective in $\Opinf$ too.
\end{proof}
For clarity, we formulate the individual statements for the two categories involved separately.
\begin{cor}\label{cor:OMOD-injectives-characterisation}
Let $A$ be a unital operator algebra and $E\in\OMODAinf$.
Then $E$ is $\Exmax$-injective in $(\OMODAinf,\Exmax)$  if and only if
$E$ is $\Exmax$-injective in $(\Opinf,\Exmax)$ and $\Exrel$-injective in $(\OMODAinf,\Exrel)$.
\end{cor}
The ``if''-part follows directly from Proposition~\ref{Prop: relative and N inj imply injective}
and the ``only if''-part follows from Proposition~\ref{Prop: injective implies others condition}
together with the CES theorem \cite[Theorem 3.3.1]{BLM04}:
for every $E\in\OMODAinf$ there exists a Hilbert space $H$ such that $B(H)\in\OMODAinf$ and $E$
is a closed submodule of~$B(H)$. The fact that $B(H)$ is the prototypical injective operator space
finished the argument.

\smallskip
This result was obtained for \C*s in \cite[Proposition 5.11]{Rosb2020}. For $E=A$ a unital \C*,
the ``if''-part is also given by \cite[Theorem~3.5]{FrankPaulsen2003}.

The next result implies \cite[Proposition~2.3]{Crann2017}, without the estimate on the constants
which are irrelevant in the completely bounded category.
It follows straight from Proposition~\ref{Prop: relative and N inj imply injective}.
\begin{cor}\label{cor:mnMOD-injectives-characterisation}
Let $A$ be a unital operator algebra and $E\in\mnMODAinf$.
Then $E$ is $\Exmax$-injective in $(\mnMODAinf,\Exmax)$ if
$E$ is $\Exmax$-injective in $(\Opinf,\Exmax)$ and $\Exrel$-injective in $(\mnMODAinf,\Exrel)$.
\end{cor}

The converse direction, however, fails in general as the following example from \cite{Crann2017} shows.
Suppose $\Gamma$ is a non-amenable discrete group; then the full \C* $C^*(\Gamma)$ is not nuclear and hence
its bidual $C^*(\Gamma)^{**}$ not injective as a \C*, hence not injective in $\Opinf$ by \cite[Theorem~3.2]{FrankPaulsen2003}.
By duality, $C^*(\Gamma)^{**}=\CB{A}{\CC}$ where $A=B(\Gamma)$ is the Fourier-Stieltjes algebra and
$\CB AG$ is always $\Exmax$-injective for an injective operator space $G$ (by the operator
algebra version of \cite[Corollary 4.11]{Rosb_arch2020}) and thus automatically $\Exrel$-injective.

Considered as a $C^*(\Gamma)$-module on the other hand, $C^*(\Gamma)^{**}$ is not injective in $\OMODinf{C^*(\Gamma)}$
by \cite[Theorem~3.4]{FrankPaulsen2003}. This reveals a subtle difference between the two operator space module categories.

\smallskip
The above results help to understand the issue of the existence of enough injectives in our module categories;
that is, for each $E\in\catsa$, do there exist $I\in\catsa$ injective and $\mu\in\Mmor(E,I)$?
When $A$ is a \C*, $B(H)$ supplies $(\OMODAinf,\Exmax)$ with enough injectives, however it loses its role when $A$
is a general operator algebra (\cite[Example~3.5]{Smith1991}) and in fact, the question seems to be open.
For \C*s, \cite[Proposition 4.13]{Rosb_arch2020} answers the question affirmatively for $(\mnMODAinf,\Exmax)$,
and the argument extends to general (unital) operator algebras.
The question remains unresolved for either of the two categories with the relative structure $\Exrel$.
The expectation seems to be that it fails, compare \cite{Aristov2000}, \cite{Aristov2002} for example.
By the above corollaries, a module which is not injective in $\Opinf$ cannot be embedded into a module which
is injective in $\Exmax$ which restricts the possible choices to modules which are only $\Exrel$-injective
and not injective in $\Exmax$ nor in $\Opinf$.

\newpage

\begin{Backmatter}

\bibliographystyle{plain}
\bibliography{refs_modules}

\printaddress

\end{Backmatter}

\end{document}